\documentclass[a4paper,12pt,draft]{article}
\usepackage{amssymb,amsmath}

\def\subsection#1{\refstepcounter{subsection}
    {\par\it\arabic{section}.\arabic{subsection} #1.}}

\begin{document}

\begin{center}
\Large On extensions and branching rules for modules close to completely splittable
\normalsize
\vskip 1.5em%
{\large Vladimir Shchigolev}\\
\medskip
vkshch@vens.ru
\end{center}

\newlength{\argwidth}

\def\eqn#1#2{\settowidth{\argwidth}{$#2$}
\hbox to \argwidth {\hss$\displaystyle#1$\hss}}

\makeatletter
\@addtoreset{equation}{section}
\renewcommand{\theequation}{\arabic{section}.\arabic{equation}}
\makeatother

\def\cs{completely splittable }
\def\acs{almost completely splittable }

\def\R{\mathop{\rm hook}\nolimits}

\def\quo{\mathop{\rm quo}}
\def\rem{\mathop{\rm rem}}

\def\#{\linebreak}
\def\sectsign{\mathhexbox278}

\def\gm{\gamma}
\def\ep{\varepsilon}
\def\f{\varphi}
\def\Gm{\Gamma}
\def\al{\alpha}
\def\bt{\beta}
\def\dl{\delta}
\def\Dl{\Delta}
\def\sg{\sigma}
\def\tt{\theta}
\def\TT{\Theta}
\def\zt{\zeta}
\def\lm{\lambda}
\def\Lm{\Lambda}
\def\om{\omega}
\def\Om{\Omega}
\def\N{\mathbb N}
\def\Z{\mathbb Z}
\def\Q{\mathbb Q}
\def\vp{\varphi}
\def\Mu{{\rm M}}
\def\Nu{{\rm N}}

\def\ch{\mathop{\rm char}}
\def\ker{\mathop{\rm Ker}}
\def\im{\mathop{\rm Im}}
\def\mod{\mathop{\rm mod}}
\def\vr{\mathop{\rm vr}}
\def\sump{\mathop{{\sum}'}}
\def\opl{\mathop{\oplus}}
\def\cf{\mathop{\rm cf}}
\def\rad{\mathop{\rm rad}}
\def\head{\mathop{\rm head}}
\def\node{\mathop{\rm node}\nolimits}
\def\Ext{\mathop{\rm Ext}\nolimits^1}
\def\Exte{\mathop{\rm Ext}\nolimits^{\rm e}}
\def\Hom{\mathop{\rm Hom}\nolimits^{\vphantom{1}}}
\def\End{\mathop{\rm End}\nolimits^{\vphantom{1}}}
\def\soc{\mathop{\rm soc}}
\def\res{\mathop{\rm res}}
\def\cont{\mathop{\rm cont}}
\def\Ind{\mathop{\rm Ind}\nolimits}
\def\Res{\mathop{\rm Res}\nolimits}
\def\Tr{\mathop{\rm Tr}\nolimits}

\def\St{\mathop{\rm St}\nolimits}
\def\st{\mathop{\rm st}\nolimits}
\def\tranc{\mathop{\rm tranc}\nolimits}

\def\b{{\rm b}}

\def\stn{\varsubsetneq}
\def\ds{\displaystyle}
\def\equ{\Leftrightarrow}

\def\<{\langle}
\def\>{\rangle}

\def\led{\vartriangleleft}
\def\ged{\vartriangleright}
\def\ledeq{\trianglelefteq}
\def\gedeq{\trianglerighteq}
\def\notled{\ntriangleleft}
\def\notged{\ntriangleright}
\def\notledeq{\ntrianglelefteq}
\def\notgedeq{\ntrianglerighteq}

\def\subsetcong{{ \lefteqn{\stackrel{\subset}{\vphantom{.}} }}_\sim}

\renewcommand{\(}{\left(}
\renewcommand{\)}{\right)}

\def\ems{\emptyset}
\def\ls{\ldots}
\def\cc{$\cdots$}
\def\upa{{\uparrow}}
\def\doa{{\downarrow}}
\def\Upa{{\dot \uparrow}}
\def\Doa{\text{\d{$\downarrow$}}}
\def\tleq{\trianglelefteq}
\def\ntleq{\ntrianglelefteq}
\def\treq{\trianglerighteq}

\def\e{\mathop{\rm e}}
\def\ind{\mathop{\rm ind}}
\def\sgn{\mathop{\rm sgn}\nolimits}
\def\s{\mathop{\rm sim}}

\def\Dok{{\bf Proof. }}
\def\enddok{\hfill$\square$}

\renewcommand{\le}{\leqslant}
\renewcommand{\ge}{\geqslant}

\makeatletter
\def\@begintheorem#1#2{\it \trivlist
       \item[\hskip \labelsep{\bf #1\ #2.}]}
\def\@opargbegintheorem#1#2#3{\it \trivlist
       \item[\hskip \labelsep{\bf #1\ #2.\ (#3)}]}
\makeatother

\newtheorem{teo}{Theorem}[section]
\newtheorem{lemma}[teo]{Lemma}
\newtheorem{gip}[teo]{Conjecture}
\newtheorem{problema}[teo]{Problem}
\newtheorem{sled}[teo]{Corollary}
\newtheorem{opr}[teo]{Definition}
\newtheorem{utv}[teo]{Proposition}

\def\codim{\qopname\relax o{codim}}

\newcounter{note}
\newcommand{\note}{\par\refstepcounter{note}%
{\bf Замечание \arabic{note}.} }
\renewcommand{\thenote}{\arabic{note}}

\makeatletter
\def\pmod#1{\allowbreak\if@display\mkern6mu\else\mkern6mu\fi({\operator@font mod}\mkern6mu#1)}
\makeatother

\renewcommand{\labelenumi}{{\rm \theenumi}}
\renewcommand{\theenumi}{(\arabic{enumi})}

\newenvironment{explicit}[1]{\par{\bf#1.} \it}{\par}

\newcounter{tabnum}
\newcommand{\tabnum}{\par\refstepcounter{tabnum}%
\medskip {\bf Таблица \arabic{tabnum}.} }
\renewcommand{\thenote}{\arabic{tabnum}}

\def\core{\mathop{\rm core}}
\def\shift{\mathop{\rm shift}}

\renewcommand\uparrow{\mathhexbox222}
\renewcommand\downarrow{\mathhexbox223}

\def\bigoplusf{{\textstyle\bigoplus}}

\def\newpi{\pi}
\renewcommand\epsilon{\varepsilon}
\def\H{\mathcal H}

\def\smallsum{\mathop\nolimits\Sigma}

\setcounter{tabnum}{0}

\renewcommand{\abstractname}{Abstract}
\begin{abstract}
We describe the modules $D^\lm\downarrow_{\Sigma_{n-1}}$ and
$D^\lm\uparrow^{\Sigma_{n+1}}$ for certain
simple $K\Sigma_n$-modules (completely splittable and some close to) $D^\lm$,
where $K$ is a field of characteristic $p>0$ and $\Sigma_n$ is
the symmetric group of degree $n$.
This result is based on an upper bound of the dimensions of
the $\Ext$-spaces between some simple modules.
\end{abstract}

\section{Introduction}\label{intro}

Fix a field $K$ of characteristic $p>0$.
We denote by $\Sigma_n$ the symmetric group of degree $n$.
We shall assume the natural inclusion $\Sigma_{n-1}\subset\Sigma_n$.
Calculation of the modules
$D^\lm\downarrow_{\Sigma_{n-1}}$ (if $n>0$) and
$D^\lm\uparrow^{\Sigma_{n+1}}$, where $D^\lm$ is the simple
$K\Sigma_n$-module corresponding to a $p$-regular partition $\lm$ of $n$,
is of great importance for the representation theory of
the symmetric group.
We have the following decomposition into blocks of
$K\Sigma_{n-1}$ and $K\Sigma_{n+1}$ (see~\cite[(55.2)]{Curtis_Reiner}
and~\cite[\sectsign1]{Kleshchev_ada}):
$$
D^\lm\downarrow_{\Sigma_{n-1}}=\bigoplus_{\alpha\in\Z_p}\Res_\alpha D^\lm,
\quad
\quad
\quad
\quad
D^\lm\uparrow^{\Sigma_{n+1}}=\bigoplus_{\alpha\in\Z_p}\Ind^\alpha D^\lm.
$$
A lot of information about
$\Res_\alpha D^\lm$ and $\Ind^\alpha D^\lm$
is contained in~\cite{Kleshchev_briv} and~\cite{Kleshchev_tf3}.
For example, the socles of these modules are known.
It is also known when a module arbitrarily chosen from these
is simple.
On the other hand, not all composition multiplicities
of nonsimple modules $\Res_\alpha D^\lm$ and $\Ind^\alpha D^\lm$
are known in the general case.

However, all the above mentioned multiplicities can be explicitly
found for some nonsimple block components.
For example, the main result of~\cite{Shchigolev11} allows us to prove
in the present paper Theorems~\ref{t:indcs:1} and~\ref{t:indcs:2}.
Together with the known Proposition~\ref{utv:prelim:1}, they yield
all composition multiplicities of
\begin{itemize}
\item $D^\lm\uparrow^{\Sigma_{n+1}}$, where
      $D^\lm$ is a completely splittable $K\Sigma_n$-module
      (see~\cite[Definition~0.1]{Kleshchev_notecs});
\item $D^{\lm^B}\downarrow_{\Sigma_n}$,
      where $D^\lm$ is a completely splittable $K\Sigma_n$-module
      distinct from $D^{(1^{p-1})}$ and $B$ is the bottom
      $\lm$-addable node.
\end{itemize}

Theorem~6 from~\cite{Shchigolev11} prompts for what other modules one
may hope to prove similar results.
A partition $\lm$ is called {\it big} (for fixed $p$),
if $D^\lm$ is completely splittable,
$\lm$ has more than one nonzero parts and
at least one rim $p$-hook.
In that case, we denote by $\tilde\lm$ the partition obtained from
$\lm$ by moving all the nodes from the last row of
the highest rim $p$-hook of $\lm$ to the first row.
Big modules exist only for $p>2$.
In the present paper, we prove an upper bound (Theorem~\ref{t:upperest:1})
of the dimensions of $\Ext_{\Sigma_n}(D^{\tilde\lm},D^\mu)$,
where $\tilde\lm$ does not strictly dominate $\mu$,
similar to the bound of~\cite[Theorem~6]{Shchigolev11}.
There are examples showing that this bound is not exact.
A separate paper is planed to prove the exact formula.
However, the bound we have obtained is enough to
\begin{itemize}
\item prove the exact formula
      (Proposition~\ref{utv:prelim:1} and Theorem~\ref{t:indacs:2})
      for composition multiplicities
      of~$D^{\tilde\lm}\uparrow^{\Sigma_{n+1}}$,
      where $\lm$ is a big partition of $n$ such that
      $h(\lm)\ge\tfrac{p+3}2$, $\lm_1+h(\lm)\ne0\pmod p$ and
      $h_{2,1}(\lm)\ne p-1$;
\item put forward conjectures on the composition multiplicities of some
      $\Ind^\alpha D^\lm$ and $\Res_\alpha D^\lm$
      (Conjectures~\ref{gip:hypothesis:1}--\ref{gip:hypothesis:3})
      confirmed by calculations
      within the known decomposition matrices.
\end{itemize}

In~\cite{Kleshchev_tf3}, the composition multiplicities of
$D^\mu$ in $\Ind^\alpha D^\lm$ and $\Res_\alpha D^\lm$
are calculated for $\mu$ one node more or respectively less
than $\lm$.
Applying the Mullineux map $m$, one can calculate the
same composition multiplicities when $m(\mu)$ is
one node more or respectively less than $m(\lm)$.
However, to prove Conjectures~\ref{gip:hypothesis:1}--\ref{gip:hypothesis:3},
one must show that the required multiplicities equal $1$
for some partitions $\mu$ having neither of the two mentioned forms.
We conjecture that this can be done by the methods of~\cite{Kleshchev_briv}.

The paper is organized as follows.
In~\sectsign2, we introduce the main objects used in the text.
The technique of the present paper differs from that of~\cite{Shchigolev11}
mainly in using abaci.
The theory of abaci is presented in~\cite{James_Kerber}.
We prove some inequalities in~\sectsign3.
In particular, the most useful inequalities
of~\cite[~\sectsign3]{Shchigolev11} are reproved
by the methods standard for homological algebra.
In~\sectsign4, an inductive method of obtaining an upper bound of
the dimensions of $\Ext_{\Sigma_n}(D^\lm,D^\mu)$ is described.
With its help, we reprove the main result of~\cite{Shchigolev11}
but in a much simpler and more visual way due to using abaci.
A slight, though necessary, modification of this method is used to
prove Theorem~\ref{t:upperest:1}.
For a more precise bound in this theorem, an auxiliary upper bound
obtained in~\sectsign7 is needed.
To sharpen this auxiliary bound, we use multiplication by the sign
representation and the Mullineux map connected with it.
The corresponding calculations are given in~\sectsign6.
Finally in~\sectsign9, we prove the above mentioned results on
the composition multiplicities of the induced and restricted modules.

\section{Notation and definitions}\label{notation}

\subsection{Generalities}\label{general}
Throughout the paper, we fix a field $K$ of positive characteristic $p$.
All rings and modules are assumed finite dimensional over $K$.
For $n\in\Z$, let $\bar n$ denote $n+p\Z$, which is an element of
$\Z_p=\Z/p\Z$.
For a pair $n,m\in\Z$, where $m>0$, let $\quo(n,m)$ and $\rem(n,m)$
denote the integers such that $n=\quo(n,m)m+\rem(n,m)$ and $0\le\rem(n,m)<m$.
For integers $r$ and $s$, the following notation will be used:
$$
\begin{array}{lr}
[r,s]=\{i\in\Z:r\le i\le s\},&(r,s]=\{i\in\Z:r<i\le s\},\\[12pt]
[r,s)=\{i\in\Z:r\le i<s\},   &(r,s)=\{i\in\Z:r<i<s\}.
\end{array}
$$

For an arbitrary assertion $\rho$, let $[\rho]$ denote $1$
if this assertion holds and denote $0$ otherwise.
For any set of integers $S$, we define its characteristic function by
$\bar S(n)=[n\in S]$ for $n\in\Z$.
We define an ordered set
$\Z'=\Z\cup\{+\infty\}$, where $+\infty>n$ holds for any $n\in\Z$.
For any module $M$ and any simple module $N$, let $[M:N]$ denote
the composition multiplicity of $N$ in $M$.

\subsection{Partitions}\label{partitions}
Given a sequence $a$, $|a|$ denotes its length.
If a positive integer $i$ is such that $i\le|a|$
in the case where $|a|<+\infty$, then $a_i$ stands for the $i^{th}$
from the beginning element of $a$.

A {\it partition} of an integer $n$ is an infinite nonincreasing sequence
of nonnegative integers, whose elementwise sum equals $n$.
To say that $\lm$ is a partition of $n$, the notation $\lm\vdash n$
is used.

In practice we write only a finite initial part of a partition
that is followed by zeros (not to be confused with finite sequences).
For example, if $\lm$ is a partition and we write $\lm=(5,3,0)$,
then $\lm_1=5$, $\lm_2=3$ and $\lm_i=0$ for $i\ge3$.
The {\it height} of a partition $\lm$ is the number $h(\lm)$ of
its nonzero entries.
A partition $\lm$ that does not contain $p$ or more identical entries is
called {\it $p$-regular}.
Let $\sum\lm$ denote the sum of all components of $\lm$.

For a partition $\lm$
we define its {\it Young diagram} by the formula
$[\lm]=\{(i,j)\in\Z\times\Z: 1\le i\le h(\lm),1\le j\le\lm_i\}$.
Elements of $\Z\times\Z$
are called {\it nodes}.
For a node $A=(i,j)$, we put $r(A)=i$
and $\res A=\overline{j-i}$.
If the diagrams of partitions $\lm$ and $\mu$ contain the same number of
nodes of each $p$-residue, then we write $\lm\sim\mu$.
{\it Removable, addable, normal, good, conormal} and {\it cogood}
nodes of a partition $\lm$ are defined in~\cite{Kleshchev_tf3}.
We also use the notations $\lm_A$ and $\lm^B$,
where $A$ is a removable and $B$ is an addable node of $\lm$,
for partitions with diagrams $[\lm]\setminus\{A\}$ and $[\lm]\cup\{B\}$
respectively.
Let $\lm^t$ denote the partition, whose diagram is
obtained by transposing $[\lm]$.

For a partition $\lm$ and an integer $i$, we put
$\sigma_i(\lm)=\sum_{1\le j\le i}\lm_j$.
A partition $\lm$ is said to dominate $\mu$
if $\sigma_i(\lm)\ge\sigma_i(\mu)$ for any $i\ge1$.
This fact is denoted by $\lm\gedeq\mu$.
The formula $\lm\ged\mu$ means that $\lm\gedeq\mu$ and $\lm\ne\mu$.

Let $h_{i,j}(\lm)$ denote the length of the hook of $\lm$
with base $(i,j)$.
We have $h_{i,j}(\lm)=\lm_i+\lm^t_j-i-j+1$.
{\it Rim, $p$-segment, $p$-edge, rim $p$-hook, $p$-core}
of a partition $\lm$ are defined in~\cite{James_Kerber} and~\cite{Mullineux1}.
Let $e(\lm)$ denote the number of nodes in the $p$-edge of a partition $\lm$
and $\vp(\lm)$ denote the partition obtained from $\lm$
by removing its $p$-edge.

\subsection{Modules}\label{modules}
To each partition $\lm$ of $n$ there corresponds a $K\Sigma_n$-module
$S^\lm$, which is called the {\it Specht module}
(see, for example, \cite[Definition~4.3]{James1}).
We put $D^\lm=S^\lm/\rad S^\lm$.
The map $\lm\mapsto D^\lm$ defines a one-to-one correspondence between
$p$-regular partitions of $n$ and simple $K\Sigma_n$-modules.

For $n\ge0$, we put $\epsilon_n=((k+1)^d,k^{p-1-d})$,
where $n=k(p-1)+d$, $k\in\Z$ and $0\le d<p-1$.
The following proposition plays an important role in the current paper.

\begin{utv}\label{u:notation:0.25}
Let $p>2$ and $\lm,\mu$ be partitions of $n$ such that
$h(\lm)<p$, $\mu$ is $p$-regular and $\lambda\not\ged\mu$.
Then we have
$$
\Ext_{\Sigma_n}(S^\lm,D^\mu)\cong
\left\{
\begin{array}{ll}
K &\mbox{if }\lm=\mu=\epsilon_n\mbox{ and }n\ge p ;\\
0  &\mbox{otherwise}.
\end{array}
\right.
$$
\end{utv}
{\bf Proof} is virtually the same as that of Theorem~2.9
from~\cite{Kleshchev_ada}.
One must only show that
$\mu=\varepsilon_n$ or $\mu=\kappa_n$
(after the assumptions $\lm_A=\epsilon_{n-1}$ and $\gamma=\epsilon_{n-1}$
are made) without using Lemma~1.6.

Indeed, $\mu_C\ledeq\gamma=\varepsilon_{n-1}=\lambda_A$ and
$\mu_C\not\led\lambda_A$, imply $\mu_C=\varepsilon_{n-1}$.
This means $\mu=\kappa_n$ or $\mu=\varepsilon_n$ or
$\mu=\varepsilon_{n-1}^D$, where $D=(p,1)$.
However, the last case is impossible, since
$\varepsilon_{n-1}^D\led\varepsilon_n,\kappa_n$ and
$\lm=\varepsilon_n$ or $\lm=\kappa_n$. \enddok

The following proposition is proved
similarly to~\cite[Theorem~2.10]{Kleshchev_ada}
but using Proposition~\ref{u:notation:0.25}
instead of~\cite[Theorem~2.9]{Kleshchev_ada}.

\begin{utv}\label{u:notation:0.5}
Let $p>2$ and $\lambda,\mu$ be partitions of $n$ such that
$h(\lm)<p$, $\mu$ is $p$-regular and $\lambda\not\ged\mu$.
Then
$\Ext_{\Sigma_n}(D^\lm,D^\mu)\cong\Hom_{\Sigma_n}(\rad S^\lm,D^\mu)$.
\end{utv}

Finally let us note the following proposition,
which follows directly from~\cite{James_Kerber}
and~\cite[Theorem~2.10]{Kleshchev_ada}.

\begin{utv}\label{u:notation:0.75}
If $\lm\not\sim\mu$ or if $p>2$, $\lm=\mu$ and $h(\lm)<p$,
then $\Ext_{\Sigma_n}(D^\lm,D^\mu)=0$.
\end{utv}

The modules $\Ind^\alpha M$ and $\Res_\alpha M$,
where $M$ is a $K\Sigma_n$-module,
are defined, for example in~\sectsign1 of~\cite{Kleshchev_ada}.

\subsection{Abaci}\label{abaci}
We shall slightly modify the classical notion of abacus
introduced in~\cite{James_Kerber} to make it more symmetrical and
convenient to work both with removing and adding nodes.
Everything what follows can be proved by the methods
of~\cite{James_Kerber} and~\cite{Kleshchev_mullpaper}.

An {\it abacus} is any map $\Lambda:\Z\to\{0,1\}$
for which there exists a number $N$ such that
$\Lm(n)=1$ for $n\le-N$ and $\Lm(n)=0$ for $n\ge N$.
The {\it shift} of an abacus $\Lm$ is the limit
$\shift(\Lm)=\lim\limits_{x\to-\infty}x+\sum_{n\ge x}\Lm(n)$
over integer $x$.
Every abacus $\Lm$ defines an injective map
$\node_\Lm:\Z\to\Z\times\Z$ by
$$
\node_\Lm(a)=\bigl(1+\sum_{n>a}\Lm(n),\sum_{n\le a}(1-\Lm(n))\bigl).
$$
We have
\begin{equation}\label{eq:notation:1}
\begin{array}{l}
c-r=a-\shift(\Lm),\text{ where }\node_\Lm(a)=(r,c);\\[6pt]
\res\node_\Lm(a)=\bar a-\overline{\shift(\Lm)}.
\end{array}
\end{equation}
Indeed, take any integer $x$ such that
$\Lm(n)=1$ for $n<x$.
Then $c-r=\sum_{x\le n\le a}(1-\Lm(n))-1-\sum_{n>a}\Lm(n)=
a-x-\sum_{x\le n}\Lm(n)=a-\shift(\Lm)$.

An element $a\in\Z$ such that $\Lambda(a)=1$ is called
a {\it bead} of $\Lambda$,
and an element $b\in\Z$ such that $\Lambda(b)=0$ is called
a {\it space} of $\Lambda$.

A bead $a$ of an abacus $\Lambda$ is called
\begin{itemize}
\itemsep=0pt
\item {\it proper} if there is a space $b$ of $\Lm$,
strictly less than $a$;
\item {\it initial} if $\Lambda(a-1)=0$;
\item {\it normal} if it is initial and
$\sum_{0<k\le s}(\Lambda(a+pk)-\Lambda(a-1+pk))\ge0$ for any $s>0$;
\item {\it good} if it is the smallest normal bead of a given $p$-residue;
\item {\it movable up} if $\Lambda(a-p)=0$.
\end{itemize}

A space $b$ of an abacus $\Lambda$ is called
\begin{itemize}
\itemsep=0pt
\item {\it initial} if $\Lambda(b-1)=1$;
\item {\it conormal} if it is initial and
$\sum_{0<k\le s}(\Lambda(b-1-pk)-\Lambda(b-pk))\ge0$ for any $s>0$;
\item {\it cogood} if it is the greatest conormal space of
a given $p$-residue.
\end{itemize}

If an abacus contains at least one proper bead,
then it is called {\it proper}.
Otherwise it is called {\it improper}.

The first formula of~(\ref{eq:notation:1}) shows
what partition should be assigned to an abacus $\Lm$.
Let $a_1$, $a_2$, \ldots be all the beads of $\Lm$
written in descending order.
We have $\node_\Lm(a_i)=(i,\lm_i)$ for some numbers $\lm_i$.
We define the partition $P(\Lm)=(\lm_1,\lm_2,\ldots)$.
$\Lm$ is said to be an {\it abacus} of $P(\Lm)$.
Let $h$ denote the number of proper beads in $\Lm$.
Clearly $P(\Lm)$ has height $h$.
In this connection the number of proper beads of an abacus
is called its {\it height}.
We have $h_{i,1}(P(\Lm))=a_i-b$ for $1\le i\le h$,
where $b$ is the smallest space of $\Lm$,
which is obviously equal to $\shift(\Lm)-h$.

Let $\Lm$ be an arbitrary abacus and $m\in\Z$.
For any $n\in\Z$, we put $(m+\Lm)(n)=\Lm(n+m)$.
An elementary verification shows that
\begin{equation}\label{eq:notation:1.5}
P(m+\Lm)=P(\Lm)\quad\mbox{ and }\quad\shift(m+\Lm)=-m+\shift(\Lm).
\end{equation}
Note that for every partition $\lm$ there is exactly one
abacus $\Lm$ of a given shift such that $\lm=P(\Lm)$.
This fact and~(\ref{eq:notation:1.5}) imply
that if $P(\Lm)=P(\Mu)$ then $\shift(\Lm)+\Lm=\shift(\Mu)+\Mu$.

\begin{utv}\label{u:notation:1}

$\node_\Lm$ bijectively maps:
\begin{enumerate}
\itemsep=0pt
\item the set of $\Lm$-initial beads to the set of $P(\Lm)$-removable nodes;
\item the set of $\Lm$-normal beads to the set of $P(\Lm)$-normal nodes;
\item the set of $\Lm$-good beads to the set of $P(\Lm)$-good nodes;
\item the set of $\Lm$-initial spaces to the set of $P(\Lm)$-addable nodes;
\item the set of $\Lm$-conormal spaces to the set of $P(\Lm)$-conormal nodes;
\item the set of $\Lm$-cogood spaces to the set of $P(\Lm)$-cogood nodes.
\end{enumerate}
\end{utv}

Let $c$ be an initial bead or an initial space of an abacus $\Lm$.
Denote by $\Lm_c$ in the former case and by $\Lm^c$ in the latter case
the abacus, whose value at $n$ is $1-\Lambda(n)$ if $n=c$ or $n=c-1$
and is $\Lm(n)$ if $n\ne c$ and $n\ne c-1$.
The operations of removing and adding nodes are connected with
the transformations of abaci just described by
$$
\begin{array}{l}
P(\Lm_c)=P(\Lm)_{\node_\Lm(c)}\mbox{ if }c\mbox{ is an initial bead},\\[10pt]
P(\Lm^c)=P(\Lm)^{\node_\Lm(c)}\mbox{ if }c\mbox{ is an initial space}.
\end{array}
$$

One of the main causes to use abaci is that they help
easily watch the removal of rim $p$-hooks and $p$-edges.
Let $a$ be a movable up bead of $\Lambda$.
Denote by $\R_\Lambda(a)$ the rim hook of $P(\Lm)$
with base $(\sum_{n\ge a}\Lm(n),\sum_{n\le a-p}(1-\Lm(n)))$.
The abacus, whose value at $n$ is $1-\Lambda(n)$
if $n=a$ or $n=a-p$ and is $\Lm(n)$ if $n\ne a$ and $n\ne a-p$,
is said to be obtained from $\Lambda$ by {\it moving $a$ one position up}.

\begin{utv}\label{u:notation:2}
The map $\R_\Lambda$ is a bijection from
the set of all movable up beads of $\Lambda$
to the set of all rim $p$-hooks of $P(\Lm)$.
Moreover, if $\bar\Lambda$ is the abacus obtained from $\Lambda$
by moving the bead $a$ one position up,
then $P(\bar\Lambda)$ is the partition obtained from $P(\Lm)$
by removing $\R_\Lambda(a)$.
\end{utv}

The above terminology is explained by
the following way of representing abaci.
Let $\Lm$ be an arbitrary abacus.
Let the position $(r,c)$ of the table $T$, where $r\in\Z$ and
$c=0$, \ldots, $p-1$, be occupied by $\cdot$ if $\Lm(pr+c)=0$ and
by $\circ$ if $\Lm(pr+c)=1$.
According to the definition of an abacus
there are two numbers $r_1$ and $r_2$ such that
row $r$ of $T$ is occupied solely by $\circ$ if $r<r_1$
and is occupied solely by $\cdot$ if $r>r_2$.
Let us draw only rows $r$ of $T$ with $r_1\le r\le r_2$
and indicate the number of any row.

In this connection the set
$\{n\in\Z:n=i\pmod p\}$, where $0\le i\le p-1$,
is called the {\it $i^{th}$ runner}.

{\it Example.} Let $\Lm=\overline{(-\infty,0]\cup\{4,8,9,11,14,17\}}$.
A representation of this abacus for $p=7$ is
$$
\begin{array}{ccccccc}
\circ&\cdot&\cdot&\cdot&\circ&\cdot&\cdot\\[-1.3pt]
\cdot&\circ&\circ&\cdot&\circ&\cdot&\cdot\\
\circ&\cdot&\cdot&\circ&\cdot&\cdot&\cdot\\
\end{array}
$$
where row $0$ is the highest.
We have $P(\Lm)=(11,9,7,6,6,3)$.

Let us continue the functions $\vp$ and $e$
defined in~\sectsign\ref{partitions}
to the set of all abaci so that
$P(\vp(\Lm))=\vp(P(\Lm))$ and $e(\Lm)=e(P(\Lm))$ for any abacus $\Lm$.

Let $\Lm$ be an arbitrary abacus.
Denote by $a$ the greatest bead and by $b$ the smallest space
of this abacus.

We put $\vp(\Lm)=\Lm$ and $e(\Lm)=0$ for any improper abacus $\Lm$.
In the remaining part of this subsection we assume that $\Lm$ is proper.
Following~\cite[Definition~1.3]{Kleshchev_mullpaper},
we consider the set of beads $\{m_1$, \ldots, $m_N\}$ of $\Lm$
(in the paper just cited it is called ``set of r-movable beads'')
that is as follows.
We put $m_1=a$.
Suppose the beads $m_1$, \ldots, $m_i$ are already chosen.
If one of the following conditions holds:
\begin{itemize}
\itemsep=0pt
\item $m_i-p$ is an improper bead;
\item $m_i-p$ is a space and there are no beads less than
      $m_i-p$,
\end{itemize}
then we put $N=i$ and we stop here.
In the opposite case, $m_{i+1}$ is equal to the greatest
bead $c$ of $\Lm$ such that $c\le m_i-p$.
It can be easily seen that all the beads $m_1$, \ldots, $m_N$
we have built are proper.

We shall build the abacus $\vp(\Lm)$ as follows.
If $m_N-p$ is an improper bead, then we move $m_N$ to position $b$.
Otherwise we move $m_N$ to position $m_N-p$.
The remaining beads $m_{N-1}$, \ldots, $m_1$ must be moved one position up
one after another starting with $m_{N-1}$ and ending with $m_1$.
Denoted by $\vp(\Lm)$ the resulting abacus.
Let
$$
e(\Lm)=\sum P(\Lm)-\sum P(\vp(\Lm)).
$$

\subsection{Operation $\H_\epsilon$}\label{Hepsilon}
In the sequel for any abacus $\Lm$ and a positive integer $i$,
$\b^\Lm(i)$ denotes the $i^{th}$ bead of $\Lm$
counting from the greatest one
and $\b_\Lm(i)$, where $i\le h(\Lm)$, denotes the $i^{th}$ proper
bead of $\Lm$ counting from the smallest one.
If $\Lm$ is proper then we abbreviate
$\b^\Lm=\b^\Lm(1)$ and $\b_\Lm=\b_\Lm(1)$.

Let $\epsilon$ be a finite set of integers.
Consider the following map
\begin{equation}\label{eq:notation:2}
\Lm+\sum_{i=1}^{|\epsilon|}
\bigl(\overline{\{\b^\Lm(i)+p\epsilon_i\}}-\overline{\{\b^\Lm(i)\}}\bigl).
\end{equation}
As we want to obtain an abacus again, we shall say
that $\H_\epsilon$ is applicable to $\Lm$ if and only
if~(\ref{eq:notation:2}) defines an abacus.
In that case we define $\H_\epsilon(\Lm)$ to be
map~(\ref{eq:notation:2}).
For a partition $\lm$, we put $\H_\epsilon(\lm)=P(\H_\epsilon(\Lm))$,
where $\lm=P(\Lm)$.
Clearly $\H_\epsilon(\lm)$ is well defined.

In other words, if $\H_\epsilon$ is applicable to $\Lm$,
then $\H_\epsilon(\Lm)$ is obtained from $\Lm$ by moving each bead
$\b^{\Lm}(i)$, where $1\le i\le|\epsilon|$,
to position $\b^\Lm(i)+p\epsilon_i$.

\begin{opr}\label{o:notation:1}
An abacus (partition) is called $\epsilon$-big,
if it is completely splittable of height $|\epsilon|$ and
$\H_\epsilon$ is applicable to it.
\end{opr}

We have the following evident assertion.

\begin{lemma}\label{l:notation:1}
Let $\Lm$ be an $\epsilon$-big abacus and $1\le i\le|\epsilon|$
such that the runner containing $\b^\Lm(i)+1$
contains no proper bead of $\Lm$.
Then $\Lm^{\b^\Lm(i)+1}$ is $\epsilon$-big and
$\H_\epsilon(\Lm^{\b^\Lm(i)+1})=\H_\epsilon(\Lm)^{\b^\Lm(i)+1+p\epsilon_i}$.
\end{lemma}

\section{Inequalities}\label{inequal}
We state the following standard fact for future reference.
\begin{utv}\label{utv:inequal:1} Let $A$, $B$, $C$ be finite dimensional
$K$-spaces and $R$ be a $K$-algebra.
\begin{enumerate}
\itemsep=0pt
\renewcommand{\labelenumi}{{\rm \theenumi}}
\renewcommand{\theenumi}{(\alph{enumi})}
\item\label{utv:inequal:1:part:1}
If $A\to B\to C$ is an exact sequence, then
$\dim B\le\dim A+\dim C$.
\item\label{utv:inequal:1:part:2}
If $A$, $B$, $C$ are $R$-modules and $A$ is isomorphic
to a submodule of $B$, then\newline $\dim\Hom_R(C,A)\le\dim\Hom_R(C,B)$.
\item\label{utv:inequal:1:part:3}
If $A$, $B$, $C$ are $R$-modules and $A$ is a homomorphic image of $B$,
then\newline $\dim\Hom_R(A,C)\le\dim\Hom_R(B,C)$.
\end{enumerate}
\end{utv}

\subsection{Case of restricted height}
We define
\begin{itemize}
\item[]
$
\varepsilon(\lm,\mu)=\left\{
\begin{array}{l}
\mbox{$1$ if $\lm=\mu^B$ for a $\mu$-conormal node $B$;}\\
\mbox{$0$ otherwise.}
\end{array}
\right.
$
\item[]
$
\gamma(\mu,\lm)=\left\{
\begin{array}{l}
\mbox{$1$ if $\mu=\lm_A$ for a $\lm$-normal node $A$;}\\
\mbox{$0$ otherwise.}
\end{array}
\right.
$
\end{itemize}

\begin{lemma}\label{l:inequal:1}
Let $p>2$,
$\lambda$ and $\mu$ be $p$-regular partitions of $n>0$,
$h(\lm)<p$, $\lm\not\ged\mu$ and $A$ be a $\mu$-good node of residue $\alpha$.
Then we have
$$
\dim\Ext_{\Sigma_n}(D^\lm,D^\mu)\le
\dim\Ext_{\Sigma_{n-1}}(\Res_\alpha D^\lm,D^{\mu_A})+
\varepsilon(\lm,\mu_A).
$$
\end{lemma}
\Dok
In view of Proposition~\ref{u:notation:0.75}, we can assume $\lm\sim\mu$.
Applying the functor $\Res_\alpha$ to the exact sequence
$0\to\rad S^\lm\to S^\lm\to D^\lm\to0$,
we get the exact sequence
$0\to\Res_\alpha\rad S^\lm\to\Res_\alpha S^\lm\to \Res_\alpha D^\lm\to0$.
Applying $\Hom_{\Sigma_{n-1}}(-,D^{\mu_A})$ to the last sequence, we
get the exact sequence
\begin{equation}\label{eq:inequal:1}
\!
\!
\!
\Hom_{\Sigma_{n-1}}(\Res_\alpha S^\lm,D^{\mu_A})\to
\Hom_{\Sigma_{n-1}}(\Res_\alpha \rad S^\lm, D^{\mu_A})\to
\Ext_{\Sigma_{n-1}}(\Res_\alpha D^\lm,D^{\mu_A}).
\!
\!
\end{equation}

By Frobenius reciprocity,~\cite[Theorem~E(iii)]{Kleshchev_tf3}
and $\lm\sim\mu$, we have
$$
\dim\Hom_{\Sigma_{n-1}}(\Res_\alpha  S^\lm ,D^{\mu_A})=
\dim\Hom_{\Sigma_n}(S^\lm,\Ind^\alpha D^{\mu_A})=\varepsilon(\lm,\mu_A).
$$
By Frobenius reciprocity, the fact that
$D^\mu$ is isomorphic to a submodule of $\Ind^\alpha D^{\mu_A}$,
Proposition~\ref{utv:inequal:1}\ref{utv:inequal:1:part:2}
and Proposition~\ref{u:notation:0.5}, we get
$$
\begin{array}{l}
\ds \dim\Hom_{\Sigma_{n-1}}(\Res_\alpha\rad S^\lm, D^{\mu_A})=
\dim\Hom_{\Sigma_n}(\rad S^\lm,\Ind^\alpha D^{\mu_A})\ge\\[10pt]
\ds \dim\Hom_{\Sigma_n}(\rad S^\lm,D^\mu)=\dim\Ext_{\Sigma_n}(D^\lm,D^\mu).
\end{array}
$$
Now it remains to apply
Proposition~\ref{utv:inequal:1}\ref{utv:inequal:1:part:1} to
sequence~(\ref{eq:inequal:1}).
\enddok

Dually we prove
\begin{lemma}\label{l:inequal:2}
Let $p>2$, $\lambda$ and $\mu$ be $p$-regular partitions of $n$,
$h(\lm)<p$, $\lm\not\ged\mu$ and $A$
be a $\mu$-cogood node of residue $\alpha$.
Then we have
$$
\dim\Ext_{\Sigma_n}(D^\lm,D^\mu)\le
\dim\Ext_{\Sigma_{n+1}}(\Ind^\alpha D^\lm,D^{\mu^A})+
\gamma(\lm,\mu^A).
$$
\end{lemma}

\subsection{Case of arbitrary height}
In the absence of the restriction $h(\lm)<p$,
the estimation of the dimension of $\Ext_{\Sigma_n}(D^\lm,D^\mu)$
can be carried out in a different way.

\begin{lemma}\label{l:inequal:3}
Let $\lambda$ and $\mu$ be distinct $p$-regular
partitions of $n>0$ and $A$ be a $\lm$-good node of residue $\alpha$.
We have
$$
\dim\Ext_{\Sigma_n}(D^\lm,D^\mu)\le
\dim\Ext_{\Sigma_{n-1}}(D^{\lm_A},\Res_\alpha D^\mu)+
\dim\Hom_{\Sigma_n}(\rad\Ind^\alpha D^{\lm_A},D^\mu).
$$
\end{lemma}
\Dok
Let $I$ be the injective hull of $D^\mu$.
Since  $\Hom_{\Sigma_n}(D^\lm,I)=0$, we have
$\Ext_{\Sigma_n}(D^\lm,D^\mu)$ $\cong\Hom_{\Sigma_n}(D^\lm,I/D^\mu)$.
Let $M$ be a $K\Sigma_n$-submodule of $I$ such that
$D^\mu\subset M$ and $M/D^\mu\cong\soc(I/D^\mu)$.
Applying the functor $\Res_\alpha$ to the exact sequence
$0\to D^\mu\to M\to M/D^\mu\to0$ and then applying
$\Hom_{\Sigma_{n-1}}(D^{\lm_A},-)$,
we get the exact sequence
\begin{equation}\label{eq:inequal:2}
\!
\!
\!
\Hom_{\Sigma_{n-1}}(D^{\lm_A},\Res_\alpha M)\to
\Hom_{\Sigma_{n-1}}(D^{\lm_A},\Res_\alpha (M/D^\mu))\to
\Ext_{\Sigma_{n-1}}(D^{\lm_A},\Res_\alpha D^\mu).
\!
\!
\end{equation}
By Frobenius reciprocity, the fact that $D^\lm$
is a homomorphic image of $\Ind^\alpha D^{\lm_A}$ and
Proposition~\ref{utv:inequal:1}\ref{utv:inequal:1:part:3}, we get
$$
\begin{array}{l}
\dim\Hom_{\Sigma_{n-1}}(D^{\lm_A},\Res_\alpha(M/D^\mu))=
\dim\Hom_{\Sigma_n}(\Ind^\alpha D^{\lm_A},M/D^\mu)\ge\\[10pt]
\dim\Hom_{\Sigma_n}(D^\lm,\soc(I/D^\mu))=
\dim\Hom_{\Sigma_n}(D^\lm,I/D^\mu)=\dim\Ext_{\Sigma_n}(D^\lm,D^\mu).
\end{array}
$$
Apply $\Hom_{\Sigma_n}(-,M)$ to the exact sequence
$$
0\to\rad\Ind^\alpha D^{\lm_A}\to\Ind^\alpha D^{\lm_A}\to D^\lm\to0
$$
and get the exact sequence
\begin{equation}\label{eq:inequal:3}
\Hom_{\Sigma_n}(D^\lm,M)\to
\Hom_{\Sigma_n}(\Ind^\alpha D^{\lm_A},M)\to
\Hom_{\Sigma_n}(\rad\Ind^\alpha D^{\lm_A},M).
\end{equation}
The first term of this sequence is $0$, since $\mu\ne\lm$

Let us see what is the image of the last morphism of
sequence~(\ref{eq:inequal:3}).
Take any $\vp\in\Hom_{\Sigma_n}(\Ind^\alpha D^{\lm_A},M)$.
Since $M/D^\mu$ is semisimple, we have
$\rad\Ind^\alpha D^{\lm_A}\subset\ker(\pi\circ\vp)$,
where $\pi:M\to M/D^\mu$ is the natural projection.
Hence $\vp(\rad\Ind^\alpha D^{\lm_A})\subset D^\mu$.

Therefore the last term of sequence~(\ref{eq:inequal:3})
can be replaced by $\Hom_{\Sigma_n}(\rad\Ind^\alpha D^{\lm_A},D^\mu)$.
Hence by Frobenius reciprocity
$$
\dim\Hom_{\Sigma_{n-1}}(D^{\lm_A},\Res_\alpha M){=}
\dim\Hom_{\Sigma_n}(\Ind^\alpha D^{\lm_A},M){\le}
\dim\Hom_{\Sigma_n}(\rad\Ind^\alpha D^{\lm_A},D^\mu).
$$
Now it remains to apply
Proposition~\ref{utv:inequal:1}\ref{utv:inequal:1:part:1}
to sequence~(\ref{eq:inequal:2}).
\enddok

Dually we prove
\begin{lemma}\label{l:inequal:4}
Let $\lambda$ and $\mu$ be distinct $p$-regular partitions of $n$
and $A$ be a $\lm$-cogood node of residue $\alpha$.
We have
$$
\dim\Ext_{\Sigma_n}(D^\lm,D^\mu)\le
\dim\Ext_{\Sigma_{n+1}}(D^{\lm^A},\Ind^\alpha D^\mu)+
\dim\Hom_{\Sigma_n}(\rad\Res_\alpha D^{\lm^A},D^\mu).
$$
\end{lemma}

{\bf Remark.} It follows from~\cite[Theorem~2.10]{Kleshchev_ada}
that Lemmas~\ref{l:inequal:3} and~\ref{l:inequal:4}
remain true for $\lm=\mu$ in the case $p>2$ and $h(\lm)<p$.

It turns out that the parameters
$\dim\Hom_{\Sigma_n}(\rad\Ind^\alpha D^{\lm_A},D^\mu)$ and\linebreak
$\dim\Hom_{\Sigma_n}(\rad\Res_\alpha D^{\lm^A},D^\mu)$
arising in Lemmas~\ref{l:inequal:3} and~\ref{l:inequal:4}
sometimes admit upper bounds.

\begin{lemma}\label{l:radind:1}
Let $\lm$ be a $p$-regular partition of $n$
and $M$ be a $K\Sigma_n$-module such that
$\Ext_{\Sigma_n}(S^\lm,\rad M)=0$, $\Hom_{\Sigma_n}(S^\lm,\rad M)=0$ and
$\head(M)\cong D^\lm$.
Then $M$ is a $K\Sigma_n$-homomorphic image of $S^\lm$.
\end{lemma}
\Dok
Applying $\Hom_{\Sigma_n}(S^\lm,-)$ to the exact sequence
$0\to\rad M\to M\to D^\lm\to 0$,
we get the exact sequence
$$
0=\Hom_{\Sigma_n}(S^\lm,\rad M)\to
\Hom_{\Sigma_n}(S^\lm,M)\to\Hom_{\Sigma_n}(S^\lm,D^\lm)\to
\Ext_{\Sigma_n}(S^\lm,\rad M)=0
$$
Hence
$\Hom_{\Sigma_n}(S^\lm,M)\cong\Hom_{\Sigma_n}(S^\lm,D^\lm)\cong K$ and
there exists a nonzero homomorphism $\vp:S^\lm\to M$.
Since $\Hom_{\Sigma_n}(S^\lm,\rad M)=0$, we have
$\im\vp\not\subset\rad M$ and thus $\im\vp=M$.
\enddok

\begin{lemma}\label{l:radind:2}
Let $\lm$ and $\mu$ be distinct $p$-regular partitions of $n>0$
such that $\Ext_{\Sigma_n}(S^\lm,D^\mu)$ $=0$.
Suppose that there exists a $\lm$-good node $A$ and a
$\mu$-good node $B$ of residue $\alpha$ and that
$A$ is the unique $\lm$-normal node of residue $\alpha$
and $\lm_A\ne\mu_B$.
Then
$\dim\Hom_{\Sigma_n}\!(\rad\Ind^\alpha D^{\lm_A},D^\mu)$ $\le\varepsilon(\lm,\mu_B)$.
\end{lemma}
\Dok
In view of Proposition~\ref{u:notation:0.75}, we can assume $\lm\sim\mu$.
We put for brevity
$n=\dim\Hom_{\Sigma_n}(\rad\Ind^\alpha D^{\lm_A},D^\mu)$.
Let $V$ be a $K\Sigma_n$-submodule of
$\rad\Ind^\alpha D^{\lm_A}$ such that
$\( \rad\Ind^\alpha D^{\lm_A} \)/V\cong\bigoplus nD^\mu$.
We put $M=\(\Ind^\alpha D^{\lm_A} \)/V$.
Applying the functor  $\Res_\alpha$ to the exact sequence
$0\to\rad M\to M\to D^\lm\to0$ and taking into account
the equivalence $\Res_\alpha D^\lm\cong D^{\lm_A}$,
which follows from the uniqueness of $A$ as a $\lm$-normal
node of residue $\alpha$, we get the exact sequence
\begin{equation}\label{eq:radind:1}
0\to\Res_\alpha\rad M\to\Res_\alpha M\stackrel{\textstyle\pi}\longrightarrow D^{\lm_A}\to0.
\end{equation}

By Frobenius reciprocity, we get
$$
0\ne\Hom_{\Sigma_n}(\Ind^\alpha D^{\lm_A},M)\cong
\Hom_{\Sigma_{n-1}}(D^{\lm_A},\Res_\alpha M).
$$
Let $\chi$ be any nonzero element of the last space.
We have
$$
\soc\Res_\alpha\rad M\cong\soc\Res_\alpha\(\bigoplusf nD^\mu\)\cong
\soc\(\bigoplusf n\Res_\alpha D^\mu\)\cong
\bigoplusf n\soc\Res_\alpha D^\mu\cong
\bigoplusf nD^{\mu_B}.
$$
It follows from this formula and from $\lm_A\ne\mu_B$
that $\im\chi\cap\Res_\alpha\rad M=0$.
Therefore $\chi$ splits the epimorphism $\pi$
of sequence~(\ref{eq:radind:1}) and
\begin{equation}\label{eq:radind:2}
\Res_\alpha M\cong D^{\lm_A}\oplus\Res_\alpha\rad M.
\end{equation}

Since $D^{\mu_B}$ is a homomorphic image of $\Res_\alpha D^\mu$,
the module $\bigoplus nD^{\mu_B}$ is a homomorphic image of $\Res_\alpha\rad M$,
which by~(\ref{eq:radind:2}) is a homomorphic image of $\Res_\alpha M$.

In view of $\head(M)\cong D^\lm$ and $\rad M\cong\bigoplus nD^\mu$,
Lemma~\ref{l:radind:1} is applicable to $M$.
Therefore $M$ is a homomorphic image of $S^\lm$.
Hence
$\Res_\alpha M$ is a homomorphic image of $\Res_\alpha S^\lm$.
As a result we get that $\bigoplus nD^{\mu_B}$
is a homomorphic image of $\Res_\alpha S^\lm$.

By Proposition~\ref{utv:inequal:1}\ref{utv:inequal:1:part:3},
Frobenius reciprocity, \cite[Theorem~E(iii)]{Kleshchev_tf3} and
$\lm\sim\mu$, we get
$$
\begin{array}{l}
\ds n=\dim\Hom_{\Sigma_{n-1}}(\bigoplusf nD^{\mu_B},D^{\mu_B})\le
\dim\Hom_{\Sigma_{n-1}}(\Res_\alpha S^\lm,D^{\mu_B})\\[10pt]
\ds=\dim\Hom_{\Sigma_n}(S^\lm,\Ind^\alpha D^{\mu_B})=
\varepsilon(\lm,\mu_B).
\end{array}
$$
\enddok

Dually we prove
\begin{lemma}\label{l:radind:3}
Let $\lm$ and $\mu$ be distinct $p$-regular partitions of $n$
such that $\Ext_{\Sigma_n}(S^\lm,D^\mu)=0$.
Suppose that there exists a $\lm$-cogood node $A$ and
$\mu$-cogood node $B$ of residue $\alpha$
and that $A$ is the unique $\lm$-conormal node of residue $\alpha$
and $\lm^A\ne\mu^B$.
Then
$\dim\Hom_{\Sigma_n}\!(\rad\Res_\alpha D^{\lm^A},D^\mu)\!$ $\le\gamma(\lm,\mu^B)$.
\end{lemma}

\subsection{Filtrations and self-duality}
For the remainder of the section, we use the notation
of~\cite[\sectsign~2]{Kleshchev_tf3}.

\begin{lemma}\label{l:spechtex:1}
Fix $\lm\in X^+(n)$ and a residue $\alpha\in\Z_p$.
Let $s_1<\cdots<s_k$ denote the set of all $j$ such that
$j$ is conormal for $\lm$ and $res (j,\lm_j+1)=\alpha$.
Take $\gamma\in X^+(n)$ such that
$$
[\Tr^\alpha L_n(\lm):L_n(\gamma)]>2\sum_{i=2}^k[\Delta_n(\lm+\varepsilon_{s_i}):L_n(\gamma)].
$$
Then $\gamma=\lm+\varepsilon_{s_1}$.
\end{lemma}
\Dok
We put $N=\Tr^\alpha L_n(\lm)$.
By~\cite[Theorem~C]{Kleshchev_tf3}, there exists
a filtration $0=N_0\subset N_1\subset\cdots\subset N_k=N$
such that for $1\le j\le k$ the module $N_i/N_{i-1}$
is a nonzero homomorphic image of $\Delta_n(\lm+\varepsilon_{s_i})$.

Suppose that $\gamma\ne\lm+\varepsilon_{s_1}$.
Then we have $[N_1:L_n(\gamma)]=[\rad N_1:L_n(\gamma)]$.
Since $N$ is contravariantly self-dual,
there exists a submodule $M\subset N$
contravariantly dual to $N/\rad N_1$.
Assumption $N_1\subset M$ (actually even $\rad N_1\subset M$)
leads to a contradiction as follows:
$$
\begin{array}{l}
\ds [N:L_n(\gamma)]=
[N/\rad N_1:L_n(\gamma)]+[\rad N_1:L_n(\gamma)]\le
2[M:L_n(\gamma)]\\[10pt]
\ds =2[N/\rad N_1:L_n(\gamma)]=2[N/N_1:L_n(\gamma)]\le
2\sum_{i=2}^k[\Delta_n(\lm+\varepsilon_{s_i}):L_n(\gamma)].
\end{array}
$$

Now the second isomorphism theorem yields
$0\ne N_1+M/M\cong N_1/N_1\cap M$, whence
$[N_1+M/M:L_n(\lm+\varepsilon_{s_1})]=1$.
Recall that $[N:L_n(\lm+\varepsilon_{s_1})]=1$ as follows
from~\cite[Theorem~B(iv)]{Kleshchev_tf3} and
$[\rad N_1:L_n(\lm+\varepsilon_{s_1})]\le
[\rad\Delta_n(\lm+\varepsilon_{s_1}):L_n(\lm+\varepsilon_{s_1})]=0$.
Now a contradiction follows from
$$
\begin{array}{l}
\ds [M:L_n(\lm+\varepsilon_{s_1})]=
[N:L_n(\lm+\varepsilon_{s_1})]-
[\rad N_1:L_n(\lm+\varepsilon_{s_1})]=1,\\[10pt]
\ds [N_1+M/M:L_n(\lm+\varepsilon_{s_1})]\le
[N:L_n(\lm+\varepsilon_{s_1})]-[M:L_n(\lm+\varepsilon_{s_1})]=0.
\end{array}
$$
\enddok

\begin{lemma}\label{l:spechtex:2.5}
Fix a $p$-regular partition $\lm\vdash r$ and a residue $\alpha\in\Z_p$.
Let $B_1,\ldots,B_k$ be all the $\lm$-conormal nodes of residue $\alpha$
counted from bottom to top.
Suppose that for some $p$-regular $\gamma\vdash r+1$ there holds
$$
[\Ind^\alpha D^\lm:D^\gamma]>2\sum_{i=2}^k[S^{\lm^{B_i}}:D^\gamma].
$$
Then $\gamma=\lm^{B_1}$.
\end{lemma}
\Dok
Choose any $n>r$.
Then~\cite[Theorem~4.16]{Kleshchev_lr10} yields
$[\Ind^\alpha D^\lm:D^\gamma]=[\Tr^\alpha L_n(\lm^t):L_n(\gamma^t)]$
and~\cite[Lemma~(6.6b)]{J.A.Green1} yields $[S^{\lm^{B_i}}:D^\gamma]=
[\Delta_n((\lm^{B_i})^t):L_n(\gamma^t)]$,
where $\lm^t$, $\gamma^t$, $(\lm^{B_i})^t$
are considered as elements of $X^+(n)$.
Now the desired result follows from Lemma~\ref{l:spechtex:1}.
\enddok

{\bf Remark.}
We conjecture that the result for $\Res_\alpha$
dual to Lemma~\ref{l:spechtex:2.5} also holds and follows
(applying~\cite[Theorems $\rm C'$, $\rm B'(iv)$]{Kleshchev_tf3})
form the lemma for $\Tr_\alpha$ reciprocal to Lemma~\ref{l:spechtex:1},
\cite[Lemma~7.4]{Kleshchev_tf3} for a suitable choice of $n$ and
\cite[Theorem 2.11(v)]{Kleshchev1}.

\section{Completely splittable partitions}\label{cs}

\subsection{General construction}\label{cs_gen_const}
Choose a set $X$, whose every element has the form
$(\lm,\mu)$, where $\lm$ and $\mu$ are $p$-regular partitions such that
$h(\lm)<p$, $\lm\not\ged\mu$, $\lm\sim\mu$ and there exists
at most one $\lm$-normal node of each residue.

Define the map $U:X\to\Z'$ inductively as follows.
We put $U(\ems,\ems)=0$.
Now let $(\lm,\mu)$ be a pair of nonempty partitions of $X$.
For each $\mu$-good node $A$, let $m_A(\lm,\mu)$ equal
the following number:
\begin{itemize}
\item $\varepsilon(\lm,\mu_A)$ if there is no $\lm$-good
      node of residue $\res A$;
\item $U(\lm_B,\mu_A)+\varepsilon(\lm,\mu_A)$
      if there is a $\lm$-good node $B$ of residue $\res A$ and
      $(\lm_B,\mu_A)\in X$.
\item $+\infty$ if there is a $\lm$-good node $B$ of residue $\res A$ and
      $(\lm_B,\mu_A)\notin X$,
\end{itemize}
Define $U(\lm,\mu)=\min\{m_A(\lm,\mu): A$ is a $\mu$-good node$\}$.

\begin{lemma}\label{l:cs:1}
Let $p>2$ and $(\lm,\mu)\in X$.
Then $\dim\Ext_{\Sigma_n}(D^\lm,D^\mu)\le U(\lm,\mu)$,
where $\lm,\mu\vdash n$.
\end{lemma}
{\bf Proof} is by induction on $n$ applying Lemma~\ref{l:inequal:1}.
\enddok

\subsection{Case of completely splittable partitions}\label{case_of_cs}

\begin{explicit}{Definition~\cite[0.1]{Kleshchev_notecs}}
An irreducible $K\Sigma_n$-module $D^\lm$ is called completely splittable
if and only if the restriction $D^\lm\doa_{\Sigma^\mu}$
to any Young subgroup $\Sigma^\mu\subset\Sigma_n$ is semisimple.
\end{explicit}
For a partition $\lm$, we put $\chi(\lm)=0$ if $\lm=\ems$ and
$\chi(\lm)=\lm_1-\lm_{h(\lm)}+h(\lm)$ otherwise.
The following result yields the exact criterion for a module to be
completely splittable.
\begin{explicit}{Theorem~\cite[2.1]{Kleshchev_notecs}}
The module $D^\lm$ is completely splittable if and only if $\chi(\lm)\le p$.
\end{explicit}

A partition $\lm$ and any abacus $\Lm$ of this partition
in the case where $D^\lm$ is completely splittable are also called
{\it completely splittable}.

The first formula of~(\ref{eq:notation:1}) shows that
for any proper abacus we have
\begin{equation}\label{eq:4:0.25}
\chi(P(\Lm))=\b^\Lm-\b_\Lm+1.
\end{equation}
It follows from this formula that any proper abacus $\Lm$
is completely splittable if and only if $\b_\Lm>\b^\Lm-p$.

\begin{opr}\label{o:cs:1}
A $(-1,0^n,1)$-big abacus (partition), where $n\ge0$,
is simply called big.
For any big abacus (partition) of height $h$, we put
$\tilde\Lm=\H_{(-1,0^{h-2},1)}(\Lm)$.
\end{opr}
It is easy to see that an abacus $\Lm$ is big if and only if
it is proper, $\b^\Lm>\b_\Lm>\b^\Lm-p$ and $\b^\Lm$ is movable up.

Let $\gamma$ be a partition and $C$ be a $\gamma$-removable node.
Then we have
\begin{equation}\label{eq:6:1}
\sigma_i(\gamma_C)=\left\{
\begin{array}{ll}
\sigma_i(\gamma)&\mbox{ if } i<r(C);\\
\sigma_i(\gamma)-1&\mbox{ if } i\ge r(C).
\end{array}
\right.
\end{equation}

\begin{lemma}\label{l:6:1}${}$\\[-18pt]
\begin{enumerate}
\itemsep=0pt
\renewcommand{\labelenumi}{{\rm \theenumi}}
\renewcommand{\theenumi}{(\alph{enumi})}
\item\label{l:6:1:c:1} If $\lm$ is a completely splittable partition,
$\mu$ is any partition, $A$ is a $\mu$-removable node,
$B$ is a $\lm$-good node,
$\res A=\res B$ and $\lm\not\gedeq\mu$, then $\lm_B\not\gedeq\mu_A$.

\item\label{l:6:1:c:2} If $\nu$ is a partition not containing
distinct removable nodes of the same residue,
$\mu$ is an arbitrary partition,
$A$ is a $\mu$-removable node, $B$ is a $\nu$-good node,
$\res A=\res B$ and $\nu\led\mu$, then $\nu_B\not\gedeq\mu_A$.
\end{enumerate}
\end{lemma}
\Dok We put for brevity $x=r(A)$ and $y=r(B)$.

\ref{l:6:1:c:1}
Suppose $\lm_B\gedeq\mu_A$.
Since $B$ is $\lm$-normal, we have $\lm_1-\lm_y+y<p$.
It follows from~(\ref{eq:6:1}) that $\sigma_i(\lm)\ge\sigma_i(\mu)$
for $i<x$ and $i\ge y$.
If $\sigma_i(\lm)\ge\sigma_i(\mu)$ also holds for each
$i$ such that $x\le i<y$, then we have $\lm\gedeq\mu$,
contrary to the hypothesis.
Therefore we assume that there exists $j$
such that $x\le j<y$, $\sigma_j(\lm)<\sigma_j(\mu)$ and
$\sigma_i(\lm)\ge\sigma_i(\mu)$ for all $i<j$.
Since $(\lm_B)_1\ge(\mu_A)_1$, we have $\lm_1\ge\mu_1-1\ge\mu_x-x$.
Hence
$$
p>\lm_1-\lm_y+y\ge(\mu_x-x)-(\lm_y-y)\ge(\mu_j-j)-(\lm_y-y).
$$
Now prove that the right hand side of
the last inequality is greater than $0$.
Suppose this is false.
The inequalities
$\sigma_j(\lm_B)\ge\sigma_j(\mu_A)$, $\sigma_j(\lm)<\sigma_j(\mu)$
and~(\ref{eq:6:1}) imply $\sigma_j(\lm)=\sigma_j(\mu)-1$ or
in a different form $\lm_j+\sigma_{j-1}(\lm)-\mu_j-\sigma_{j-1}(\mu)+1=0$.
Hence
$$
\bigl(\sigma_{j-1}(\lm)-\sigma_{j-1}(\mu)\bigl)+\bigl(\lm_j-\lm_y+1\bigl)+
\bigl(y-j\bigl)\le0.
$$
This formula is a contradiction,
as the expression in its first pair of brackets is nonnegative and
the expressions in the other two are positive.
Therefore, we have $p>(\mu_x-x)-(\lm_y-y)>0$ and $\res A\ne\res B$.

\ref{l:6:1:c:2}
Suppose $\nu_B\gedeq\mu_A$.
It follows from~(\ref{eq:6:1}) that
$\sigma_i(\nu)\ge\sigma_i(\mu)$ for $i<x$ and $i\ge y$.
If $\sigma_i(\nu)\ge\sigma_i(\mu)$ also holds
for each $i$ such that $x\le i<y$,
then we have $\nu\gedeq\mu$, contrary to the condition $\nu\led\mu$.
Therefore we assume that there exists $j$ such that
$x\le j<y$, $\sigma_j(\nu)<\sigma_j(\mu)$ and
$\sigma_i(\nu)\ge\sigma_i(\mu)$ for all $i<j$.
Since $\nu\led\mu$, we have $\sigma_i(\nu)=\sigma_i(\mu)$ for all $i<j$
and thus $\nu_i=\mu_i$ for all $0\le i<j$.
Similarly to part~\ref{l:6:1:c:1}, we get $\sigma_j(\nu)=\sigma_j(\mu)-1$,
whence $\nu_j=\mu_j-1$.

If $x<j$ then $A$ and $B$ are distinct $\nu$-removable nodes
of the same residue, which is a contradiction.
Therefore $x=j$ and $A$ is $\nu$-addable.
By the definition of a normal node there exists a $\nu$-removable node
of residue $\res B$ strictly between $B$ and $A$.
This is again a contradiction.
\enddok

Let us apply the construction described in~\sectsign\ref{cs_gen_const}
to the following set
$$
X=\{(\lm,\mu):\lm\mbox{ is completely splittable},\;%
              \mu\mbox{ is $p$-regular},\;%
              \lm\not\gedeq\mu,\;%
              \lm\sim\mu\},%
$$
which we fix until the end of this section.
Note that $X=\ems$ if $p=2$, as
in that case the following conditions cannot hold
simultaneously:
$\lm$ is completely splittable, $\lm\not\gedeq\mu$, $\lm\sim\mu$.

\begin{lemma}\label{l:cs:2}
Let $(\lm,\mu)\in X$.
Then $U(\lm,\mu)=0$ except the case $\mu=\tilde\lm$,
in which we have $U(\lm,\tilde\lm)\le1$.
\end{lemma}
\Dok
Induction on $n=\sum\lm$.
For $n=0$, by definition we have $U(\lm,\mu)=0$.

Now let $n>0$.
Suppose that the theorem is true for partitions of numbers less than $n$.
Choose some abaci $\Mu$ and $\Lm$ of the same shift
such that $\mu=P(\Mu)$ and $\lm=P(\Lm)$.
Let $A$ be an arbitrary $\mu$-good node.
It exists because $\mu$ is nonempty and $p$-regular.
Denote by $a$ the bead of $\Mu$ corresponding to $A$.

First consider the case $\epsilon(\lm,\mu_A)=1$.
Then $\Lm=(\Mu_a)^c$ for some conormal space $c$ of $\Mu_a$.
Since $\lm\sim\mu$ and $\lm\not\gedeq\mu$, $a$ and $c$ are in the same runner
and moreover $a$ is below $c$.
We have $\Mu=(\Lm_c)^a$.
This at once implies $c=\b_\Lm$, $a=\b_\Lm+p$, $\b^\Lm+1=\b_\Lm+p$ (i.e. $\chi(\lm)=p$)
and $\Mu=\tilde\Lm$ (i.e. $\mu=\tilde\lm$).
Note that in this case the only bead of $\Lm$ from the same runner as $a$
is $c$, which is not normal.
Therefore, there is no $\lm$-good node of residue $\res A$
and $U(\lm,\mu)\le m_A(\lm,\mu)=\epsilon(\lm,\mu_A)=1$.

Now consider the case $\epsilon(\lm,\mu_A)=0$
If there is no $\lm$-good node of residue $\res A$,
then $U(\lm,\mu)=m_A(\lm,\mu)=0$.
Therefore we assume that there is a $\lm$-good node $B$
of residue $\res A$.
Applying Lemma~\ref{l:6:1}\ref{l:6:1:c:1},
we get $\lm_B\not\gedeq\mu_A$, whence $(\lm_B,\mu_A)\in X$.
Thus $U(\lm,\mu)\le m_A(\lm,\mu)=U(\lm_B,\mu_A)$.
Therefore we shall consider the case $U(\lm_B,\mu_A)>0$.

Let $b$ be the bead of $\Lm$ corresponding to $B$.
By the inductive hypothesis, we get $\mu_A=\widetilde{\lm_B}$,
$\Mu_a=\widetilde{\Lm_b}$ and $U(\lm_B,\mu_A)\le1$.
Hence $s\le\b^{\Lm_b}-p$, where $s$ is the smallest space of $\Lm_b$.
Therefore $b\ne s$, as otherwise we would
get $\b^\Lm-\b_\Lm=\b^{\Lm_b}-s\ge p$,
contrary to the fact that $\Lm$ is completely splittable.

Recall that $a$ and $b$ are cogood spaces of $\Mu_a$ and $\Lm_b$
respectively,
since $a$ and $b$ are good beads of $\Mu$ and $\Lm$ respectively.
We have $b=\b^{\Lm_b}(i)+1<\b_{\Lm_b}+p$ for $1\le i\le h(\Lm)$,
since $\b_{\Lm_b}+p$ is not a cogood space of $\Lm_b$.
Now Lemma~\ref{l:notation:1} can be applied to $\Lm_b$.
This yields $\tilde\Lm=\Mu_a^{b+p\epsilon_i}$, where
$\epsilon=(-1,0^{h(\Lm)-2},1)$.
Since $a$ and $b$ are in the same runner and
the only cogood space of $\Mu_a$, which is equal to $\widetilde{\Lm_b}$,
being in the same runner as $b$ is the space $b+1+p\epsilon_i$,
then $a=b+1+p\epsilon_i$ and $M=\tilde\Lm$.
\enddok

\subsection{Exact formula}\label{exact_formula_for_cs}
The next lemma shows that Definition~\ref{o:cs:1} of a big partition
and the map $\lm\mapsto\tilde\lm$ given in the current paper
is equivalent to~\cite[Definition~4]{Shchigolev11}.

\begin{lemma}\label{l:cs:3}
A partition $\lm$ is big if and only if
$\lm$ is completely splittable of height more than one
and $h_{1,1}(\lm)\ge p$.
In that case $[\tilde\lm]$ is obtained from $[\lm]$ by moving all the nodes
of the rim $p$-hook with base $(1,\lm_1+h(\lm)-p)$
from the last row to the first row.
\end{lemma}
\Dok
The first part follows from $h_{1,1}(\lm)=\b^\Lm-c$,
where $\Lm$ is a proper abacus and $c$ is its smallest space.

Now let $\lm$ be a big partition.
We put for brevity $h=h(\lm)$ and $j=\lm_1+h-p$.
Since $1\le j\le\lm_1$, the node $(1,j)$ belongs to $[\lm]$.
We have $\lm^t_j\le h$.
The inequality $\lm^t_j<h$ would imply $\lm_h<j$ and $\chi(\lm)>p$.
Hence $\lm^t_j=h$ and $h_{1,j}(\lm)=\lm_1-j+\lm^t_j=p$.

Since $\b^\Lm-p<\b_\Lm$ and $\b^\Lm<\b_\Lm+p$,
for each bead $d$ of $\Lm$ such that $\b_\Lm<d<\b^\Lm$,
the numbers of beads preceding $d$ in $\Lm$ and $\tilde\Lm$
coincide and the numbers of beads following $d$ in $\Lm$ and $\tilde\Lm$
also coincide.
The beads $\b^\Lm-p$ and $\b_\Lm+p$ are respectively the smallest
and the greatest beads of $\tilde\Lm$.
There are $\b^\Lm-p-c=\lm_1-1+h-p$ spaces in $\tilde\Lm$
preceding the former bead and
$\b_\Lm+p-c-(h-1)=\lm_h-h+p+1$ spaces
preceding the latter bead
($c$ is the smallest space of $\Lm$).
Thus $\tilde\lm=(\lm_h-h+p+1,\lm_2$, \ldots, $\lm_{h-1},\lm_1+h-p-1)$.
This is exactly the partition, whose diagram is obtained from $[\lm]$
by moving all the nodes of the rim $p$-hook with base
$(1,j)$ from the last row to the first row.
\enddok

\begin{teo}[\mbox{ \cite[Theorem~6]{Shchigolev11} }]\label{t:cs:1}
Let $p>2$ and $\lm,\mu$ be $p$-regular partitions of $n$ such that
$D^\lm$ is completely splittable and $\lm\not\ged\mu$.
Then
$$
\Ext_{\Sigma_n}(D^\lambda,D^\mu)\cong
\left\{
\begin{array}{ll}
K &\mbox{if }\mu=\tilde\lambda;\\
0 &\mbox{otherwise}.
\end{array}
\right.
$$
\end{teo}
\Dok
The equality $\Ext_{\Sigma_n}(D^\lambda,D^\mu)=0$ in the case where
$\lm$ is not big or $\mu\ne\tilde\lm$
follows from Lemma~\ref{l:cs:2} and Proposition~\ref{u:notation:0.75}.

Now suppose $\mu=\tilde\lm$.
Lemma~\ref{l:cs:3} implies that the length of the hook $h_{i,j}(\lm)$,
where $j=\lm_1+h(\lm)-p$ and $1\le i\le\lm^t_j=h(\lm)$
is divisible by $p$ if and only if $i=1$.
Since $h(\lm)>1$, the Carter conjecture proved in~\cite{James3} implies that
the module $S^\lm$ is not simple.
By Lemma~\ref{l:cs:2}, a unique top composition factor
of a nonzero module $\rad S^\lm$ is $D^{\tilde\lm}$.
This fact and the second assertion of Lemma~\ref{l:cs:2} yield
$\dim\Ext_{\Sigma_n}(D^\lm,D^{\tilde\lm})=1$.
\enddok

\subsection{Almost completely splittable partitions}\label{acs}

\begin{opr}\label{o:acs:1}
An abacus (partition) $\Mu$ is called almost completely splittable
if $\Mu=\tilde\Lm$ for some completely splittable $\Lm$,
which is called the preimage of $\Mu$.
\end{opr}
It is important to notice that the preimage $\Lm$
is uniquely determined by $\Mu$.
Indeed, let $a=\b^\Mu$.
If $\Mu$ is a completely splittable abacus, then denote by $b$ its greatest improper bead.
If $\Mu$ is not a completely splittable abacus, then put $b=\b_\Mu$.
Then $\Lm$ is obtained from $\Mu$ by moving $a$ up and moving $b$ down
one position.

A module $D^\lm$, where $\lm$ is an almost completely splittable partition,
is also called {\it almost completely splittable}.

It is convenient to define many abaci encountered in this paper
with the help of the following construction.
Take $x\in\Z$ and $S\subset{[0,p-1]}$.
Let $x_0$, \ldots, $x_i$, \ldots be all the elements of
the set $\{n\in\Z:n\ge x, \rem(n,p)\in S\}$
written in the ascending order.
For $i\ge0$, we put $\<x,S,i\>=\{x_k:i\le k<i+|S|\}$.
Note that if $i>0$ then the set $\<x,S,i-1\>$
is obtained from $\<x,S,i\>$ by replacing the greatest element $a$
of the latter set with $a-p$ (moving $a$ one position up).

{\it Example.} Let $p=7$ and $S=\{1,3,4,6\}$.
We have $\<9,S,5\>=\{18,20,22,24\}$ and $\<9,S,6\>=\{20,22,24,25\}$.

\begin{opr}
Let $h$, $i$, $x$ be integers such that
$0<i\le h<p$ and $x\ge0$.
Define
$\ds
\Lm^{(h,i,x)}=\overline{(-\infty,0)\cup\<0,[0,i),x\>\cup[i,h)}
$
and
$
\ds\lm^{(h,i,x)}=P(\Lm^{(h,i,x)}).
$
\end{opr}

Clearly, the partition $\lm^{(h,i,x)}$ has $p$-weight $x$
and therefore is a partition of $px$.
Moreover, if $x>0$ the height of $\lm^{(h,i,x)}$ is $h$,
$\core(\lm^{(h,i,x)})=\ems$, the top removable node of
$\lm^{(h,i,x)}$ is its only normal node and this node has residue $-\bar h$.
Obviously the equalities $\lm^{(h,i,0)}=\ems$ and
$\lm^{(h,i,x)}=\lm^{(h,x,x)}$
(following from $\Lm^{(h,i,x)}=\Lm^{(h,x,x)}$) for $x<i$
represent all nontrivial equalities between partitions $\lm^{(h,i,x)}$.
We put for brevity $\lm^{(H,x)}=\lm^{(H,H,x)}$.
Explicitly $\Lm^{(h,i,x)}$ is written as
$$
\Lm^{(h,i,x)}=\overline{(-\infty,0)\cup[(q+1)p,(q+1)p+r)\cup[qp+r,qp+i)\cup[i,h)},
$$
where $q=\quo(x,i)$ and $r=\rem(x,i)$.

{\bf Remark}. In~\cite{Hemmer1},
the partitions $\lm^{(H,x)}$ are called minimal.

\begin{utv}\label{utv:3:1}
If $\lm$ is a completely splittable partition, then $\chi(\core(\lm))<p$.
\end{utv}
{\bf Proof} follows from~(\ref{eq:4:0.25})
and the fact that $\core(\lm)$ is completely splittable.
\enddok

\begin{lemma}\label{l:3:1}
Let $\lm$ be a completely splittable partition, $\chi(\lm)=p$ and
the residue of any $\core(\lm)$-normal node
is equal to the residue of the bottom $\lm$-removable node.
Then $\lm=\lm^{(H,x)}$, where $1<H<p$ and $H\nmid x$.
\end{lemma}
\Dok
Suppose that $\core(\lm)$ is not empty.
Proposition~\ref{utv:3:1} implies that
$\chi(\core(\lm))<p$ and $\core(\lm)$ contains
only one removable node.
Hence $\core(\lm)=P\bigl(\overline{(-\infty,0)\cup[b,a]}\bigl)$, where $0<b\le a<p$.
Therefore
$\lm=P\bigl(\overline{(-\infty,0)\cup[(q+1)p+b,(q+1)p+r)\cup[qp+r,qp+a]}\bigl)$
for $q$ and $r$ such that $q\ge0$ and $b\le r\le a$.

The case $r=b$ is impossible, as we would get
the contradiction $\chi(\lm)=\chi(\core(\lm))<p$.
The case $r>b$ is also impossible.
Indeed, in this case the smallest initial bead of
$\overline{(-\infty,0)\cup[(q+1)p+b,(q+1)p+r)\cup[qp+r,qp+a]}$,
which is equal to $qp+r$, and the normal bead $b$ of
$\overline{(-\infty,0)\cup[b,a]}$ would belong to different runners.
This is a contradiction.

Now the desired assertion follows from
$\core(\lm)=\emptyset$ and $\chi(\lm)=p$.
\enddok

\begin{utv}\label{utv:3:2}
Let $\lm$ be a nonempty partition and $\Lm$ be its abacus.
Then $\lm-(1^{h(\lm)})=P(\Lm')$, where $\Lm'$ is obtained from $\Lm$
by replacing its smallest space by the bead.
\end{utv}

\begin{sled}\label{sl:3:1}
Let $\lm$ and $\mu$ be partitions such that
$h(\lm)=h(\mu)$ and $\core(\lm)=\core(\mu)$.
Then $\core(\lm-(1^{h(\lm)}))=\core(\mu-(1^{h(\mu)}))$.
\end{sled}

\begin{lemma}\label{l:3:2}
Let $\chi(\lm)=p$, $\core(\lm)=\core(\mu)=\ems$, $h(\mu)\le h(\lm)<p$
and the residue of any $\mu$-normal node
is equal to the residue of the bottom $\lm$-removable node.
Then $h(\mu)<h(\lm)$.
\end{lemma}
\Dok
Suppose $h(\mu)=h(\lm)$.
We put $h=h(\lm)$ and $n=\sum\lm$.

Case $n>p$.
We have $\lm=P(\Lm^{(h,x)})$.
We put $\alpha=\res(h,\lm_h)$, $\bar\lm=\lm-(1^h)$ and $\bar\mu=\mu-(1^h)$.
In the case under consideration $x>1$.
Therefore $\lm_h>1$ and $\chi(\bar\lm)=p$.
By Corollary~\ref{sl:2:0}, for any $\core(\bar\mu)$-normal node, there is
a $\bar\mu$-normal node of the same residue.
However all $\bar\mu$-normal nodes have residue $\alpha-1$,
which is equal to the residue of the bottom $\bar\lm$-removable
node $(h,\lm_h-1)$.
Since $\core(\bar\mu)=\core(\bar\lm)$
by Corollary~\ref{sl:3:1}, the residue of any $\core(\bar\lm)$-normal node
is equal to the residue of the bottom $\bar\lm$-removable node.
By Lemma~\ref{l:3:1}, we get $\core(\bar\lm)=\ems$, which is a contradiction,
since $\bar\lm$ is a partition of the number $n-h$ not divisible by $p$.

Case $n=p$ is reduced to the previous one by considering the pair
of partitions \#$\hat\lm=\lm+(p^h)$ and $\hat\mu=\mu+(p^h)$.
\enddok

\section{Removal of locally highest $p$-hooks}\label{hookrem}

\begin{opr}\label{o:2:3}
Let $u$ be a movable up bead of $\Lambda$.
The rim $p$-hook corresponding to it is called locally highest if
$u+1$ is not a movable up bead.
\end{opr}

\begin{lemma}\label{l:2:1}
Let $\lm$ be a partition and $\bar\lm$ be the partition obtained from $\lm$
by removing its locally highest rim $p$-hook $R$.
Then for any $\alpha\in\Z_p$ the number of
$\bar\lm$-normal nodes of residue $\alpha$ does not exceed
the number of $\lm$-normal nodes of the same residue.
\end{lemma}
\Dok
Let $\Lambda$ and $\bar\Lambda$ be some abaci of $\lm$ and $\bar\lm$
respectively having shifts of $p$\,-residue $1-\alpha$.
Consider the decomposition $\R_\Lambda^{-1}(R)=pi+j$, where $0\le j<p$.

By Proposition~\ref{u:notation:1} and
the second formula of~(\ref{eq:notation:1}), the theorem
will be proved if we define an embedding $\iota$
of the set of normal beads of the first runner of $\bar\Lambda$ into
the set of normal beads of the first runner of $\Lambda$.
If $j>1$ then we can take the identity map for $\iota$.
Thus we assume $0\le j\le 1$.

First consider the case $j=0$.
Let $a=px+1$ be a normal bead of $\bar\Lambda$.

If $x>i$ then $a$ is obviously a normal bead of $\Lambda$.
In that case we put $\iota(a)=a$.

If $x=i$ then $a-p$ is a bead of $\Lambda$,
as otherwise it would be possible to move the bead $a$ one
position up in $\Lambda$, which contradicts the fact that
the rim $p$-hook we have removed is locally highest.
Take any $s>i-1$.
We have
$$
\sum_{i-1<k\le s}
(\Lambda(pk+1)-\Lambda(pk))=\sum_{i<k\le s}(\Lambda(pk+1)-\Lambda(pk))=
\sum_{i<k\le s}(\bar\Lambda(pk+1)-\bar\Lambda(pk))\ge0.
$$
Now it is clear that $a-p$ is a normal bead of $\Lambda$.
We put $\iota(a)=a-p$.

The case $x=i-1$ is impossible, as otherwise
the bead $a$ would not be initial in $\bar\Lambda$.

Finally let $x<i-1$.
For any $s>x$, we have
$\sum_{x<k\le s}\bar\Lambda(pk)\ge\sum_{x<k\le s}\Lambda(pk)$ and
therefore
$$
\sum_{x<k\le s}(\Lambda(pk+1)-\Lambda(pk))\ge\sum_{x<k\le s}(\bar\Lambda(pk+1)-\bar\Lambda(pk))\ge0.
$$
Thus $a$ is a normal bead of $\Lambda$.
We put $\iota(a)=a$.

Now consider the case $j=1$.
Define the parameter $x_0$ as follows.
Let $px+1$ and $py+1$ be normal beads of $\bar\Lambda$
such that $y<x<i-1$.
We have
\begin{equation}\label{eq:2:1}
\begin{array}{l}
\ds\sum_{y<k\le i-1}(\bar\Lambda(pk+1)-\bar\Lambda(pk))=\\[20pt]
\ds\sum_{y<k\le x-1}(\bar\Lambda(pk+1)-\bar\Lambda(pk))+1+
\sum_{x<k\le i-1}(\bar\Lambda(pk+1)-\bar\Lambda(pk))\ge1.
\end{array}
\end{equation}
Therefore there is at most one number $x_0<i-1$ such that
$px_0+1$ is a normal bead of $\bar\Lambda$ and
\begin{equation}\label{eq:2:2}
\sum_{x_0<k\le i-1}(\bar\Lambda(pk+1)-\bar\Lambda(pk))=0.
\end{equation}
If there is no such number at all, then we put $x_0=+\infty$.
Under this definition, we get

\begin{explicit}{Basic property of $x_0$}
Let $px+1$ be a normal bead of $\bar\Lambda$ such that $x<i-1$.
Then $x\le x_0$ and
$\sum_{x<k\le i-1}(\bar\Lambda(pk+1)-\bar\Lambda(pk))>0$ if $x\ne x_0$.
\end{explicit}

Indeed in the case $x_0<+\infty$ this fact follows from~(\ref{eq:2:1}).
In the case $x_0=+\infty$ the sum under consideration
is not equal to zero by the definition of $x_0$
and thus is strictly positive.

If $x>i$ then similarly to the previous case we get that
$a$ is a normal bead of $\Lambda$.
We put $\iota(a)=a$.

The case $x=i$ is impossible, as $\bar\Lambda(pi+1)=0$ and
$\bar\Lambda(px+1)=1$.

If $x=i-1$ then $\Lambda(pi)=\bar\Lambda(pi)=0$, because
$\bar\Lambda(pi+1)=0$ and $a=p(i-1)+1$ is a normal bead of $\bar\Lambda$.
Take any $s>i$.
We have
$$
\sum_{i<k\le s}(\Lambda(pk+1)-\Lambda(pk))=\sum_{i-1<k\le s}(\bar\Lambda(pk+1)-\bar\Lambda(pk))\ge0.
$$
Now it is clear that $pi+1$ is a normal bead of $\Lambda$.
We put $\iota(a)=pi+1$.

If $x=x_0$ then by property~(\ref{eq:2:2}) we have
$$
\begin{array}{l}
\ds 0\le\sum_{x_0<k\le i}(\bar\Lambda(pk+1)-\bar\Lambda(pk))=-\Lambda(pi)\\[20pt]
\ds 0=\sum_{x_0<k\le i-1}(\bar\Lambda(pk+1)-\bar\Lambda(pk))=\\[20pt]
\ds \sum_{x_0<k\le i-2}(\bar \Lambda(pk+1)-\bar\Lambda(pk))+1-\Lambda(p(i-1))\ge1-\Lambda(p(i-1)).
\end{array}
$$
Therefore $\Lambda(pi)=0$ and $\Lambda(p(i-1))=1$.

Take any $s>i$.
We have
$$
\sum_{i<k\le s}(\Lambda(pk+1)-\Lambda(pk))=
\sum_{x_0<k\le s}(\bar\Lambda(pk+1)-\bar\Lambda(pk))\ge0.
$$
For the last rearrangement, we used property~(\ref{eq:2:2}) and the equalities
$\bar\Lambda(pi)=\bar\Lambda(pi+1)=0$.
Therefore $pi+1$ is a normal bead of $\Lambda$.
We put $\iota(a)=pi+1$.

Finally let $x<i-1$ and $x\ne x_0$.
Clearly, for an arbitrary $s>x$ not equal to $i-1$ we have
$$
\sum_{x<k\le s}(\Lambda(pk+1)-\Lambda(pk))=
\sum_{i<k\le s}(\bar\Lambda(pk+1)-\bar\Lambda(pk))\ge0.
$$
On the other hand by the basic property of $x_0$ we have
$$
\sum_{x<k\le i-1}(\Lambda(pk+1)-\Lambda(pk))=
\sum_{x<k\le i-1}(\bar\Lambda(pk+1)-\bar\Lambda(pk))-1\ge0.
$$
Therefore $a$ is a normal bead of $\Lm$.
We put $\iota(a)=a$.

The preimage of a bead $b$ belonging to the image of $\iota$,
is given by the following formulas:\\[12pt]
\begin{tabular}{ll}
Case $j>1$:&$\iota^{-1}(b)=b$;\\[12pt]
Case $j=0$:&
$\iota^{-1}(b)=
\left\{
\begin{array}{ll}
b    & \mbox{ if }b\ne p(i-1)+1;\\
pi+1& \mbox{ if }b=p(i-1)+1,\\
\end{array}
\right.$
\\[20pt]
Case $j=1$:&
$\iota^{-1}(b)=
\left\{
\begin{array}{ll}
b      & \mbox{ if }b\ne pi+1;\\
p(i-1)+1& \mbox{ if }b=pi+1\mbox{ and }\Lambda(p(i-1))=0;\\
px_0+1& \mbox{ if }b=pi+1\mbox{ and }\Lambda(p(i-1))=1.\hspace{35pt}\square\\
\end{array}
\right.$
\end{tabular}

\begin{sled}\label{sl:2:0}
Let $\lm$ be a partition.
For any $\alpha\in\Z_p$, the number of
$\core(\lm)$-normal nodes of residue $\alpha$
does not exceed the number of $\lm$-normal nodes of the same residue.
\end{sled}
{\bf Proof} follows from the fact that $\core(\lm)$
can be obtained by removing the highest rim $p$-hook at each step. \enddok

\section{Mullineux map of some partitions}\label{mullmap}

To calculate the Mullineux map of a partition $\lm$,
we shall use the Mullineux symbol defined in~\cite{Bessenrodt_Olsson},
which is the array
$$
G_p(\lm)=\binom{A_0\cdots A_z}{R_0\cdots R_z},
$$
where $A_j=e(\vp^j(\lm))$, $R_j=h(\vp^j(\lm))$ and $\vp^{z+1}(\lm)=\ems$.
The product of such arrays is understood as follows:
$$
\binom{A_0\cdots A_z}{R_0\cdots R_z}\binom{A'_0\cdots A'_{z'}}{R'_0\cdots R_{z'}}=
\binom{A_0\cdots A_zA'_0\cdots A'_{z'}}{R_0\cdots R_zR'_0\cdots R_{z'}}.
$$

\begin{lemma}\label{l:4:2}
Let $1<H<p$, $x>0$ and $H\nmid x$.
We put $Q=\quo(x,H)$ and $R=\rem(x,H)$.
If $R>1$ then

\medskip

$\ds
\quad
\quad
\quad
G_p\left(\widetilde{\lm^{(H,x)}}\right)=
\binom{2p   }{H}^Q
\binom{p+R-1}{H}
\binom{p    }{H-1}^{x-2(Q+1)}
\binom{p-R+1}{H-R+1}
$

\medskip

\noindent and if $R=1$ then

\medskip

$\ds
\quad
\quad
\quad
G_p\left(\widetilde{\lm^{(H,x)}}\right)=
\binom{2p}{H}^Q
\binom{p}{H-1}^{x-2Q}.
$
\end{lemma}
\Dok
Let $S=[0,H-1]\setminus\{R-1\}$,
$\Lm=\widetilde{\Lm^{(H,x)}}$ and $\lm=\widetilde{\lm^{(H,x)}}$.
It is easy to check that
$$
\vp^j(\Lm)=
\overline{(-\infty,0)\cup\<0,\{R-1\},Q-j\>\cup\<0,S,x-Q-j\>}
\quad
\mbox{ for }
\quad
1\le j\le Q.
$$
Hence $h(\vp^j(\lm))=H$
and $e(\vp^j(\lm))=2p$ for $0\le j<Q$.

Case $R>1$.
Then $h(\vp^Q(\lm))=H$.
We have
$$
\vp^{Q+j}(\Lm)=\overline{(-\infty,0]\cup\<1,S,x-2Q-1-j\>}
\quad
\mbox{ for }
\quad
1\le j\le x-2Q-1.
$$
Hence $e(\vp^Q(\lm))=p+R-1$,
$h(\vp^{Q+j}(\lm^{(H,x)}))=H-1$, $e(\vp^{Q+j}(\lm^{(H,x)}))=p$
for $1\le j<x-2Q-1$ and
$h(\vp^{x-Q-1}(\lm^{(H,x)}))=H-R+1$.
Finally $\vp^{x-Q}(\lm^{(H,x)})=\ems$ and
$e(\vp^{x-Q-1}(\lm^{(H,x)}))=p-R+1$.

Case $R=1$.
Then $h(\vp^Q(\lm))=H-1$.
We have
$$
\vp^{Q+j}(\Lm)=\overline{(-\infty,0]\cup\<1,S,x-2Q-j\>}
\quad
\mbox{ for }
\quad
0\le j\le x-2Q.
$$
Hence $h(\vp^{Q+j}(\lm))=H-1$ for $1\le j<x-2Q$,
$\vp^{x-Q}(\lm)=\ems$ and
$e(\vp^{Q+j}(\lm))=p$ for $0\le j<x-2Q$.
\enddok

\begin{opr}\label{o:4:2}
Let $H$ and $x$ be integers such that $1<H<p$, $x>0$ and $H\nmid x$.
Define
$$
\begin{array}{lcl}
\Nu^{(H,x)}&=&
\overline{(-\infty,H)\cup\<H,S_1,Q\>
\cup\<p+H-R,S_2,x-Q\>},\\[12pt]
\nu^{(H,x)}&=&P(\Nu^{(H,x)}),
\end{array}
$$
where $Q=\quo(x,H)$, $R=\rem(x,H)$, $S_1=[0,p-1]\setminus\{H-R\}$ and
$S_2=\{H-R\}\cup [H,p-1]$.
\end{opr}

To gain a better understanding of the structure of $\Nu^{(H,x)}$,
we introduce the following notation:
let $a^ {(H,x)}_0$, $a^{(H,x)}_1$, \ldots be all the elements of the set
$
\{n\in\Z:n\ge H, \rem(n,p)\in S_1\}
$
and $b^{(H,x)}_0$, $b^{(H,x)}_1$, \ldots be all the elements of the set
$
\{n\in\Z:n\ge p+H-R, \rem(n,p)\in S_2\}
$
written in ascending order.
An easy verification shows that
\begin{equation}\label{eq:4:1}
\begin{array}{l}
\ds a^{(H,x)}_y=H+y+\left[\tfrac{R-1+y}{p-1}\right],\\[16pt]
\ds b^{(H,x)}_y=p+H-1+y+(H-R)\left[\tfrac{y}{p-H+1}\right]+
(R-1)\left[\tfrac{y-1}{p-H+1}\right].
\end{array}
\end{equation}

\begin{lemma}\label{l:4:3}
Let $H$ and $x$ be integers such that
$1<H<p$, $x>0$ and $H\nmid x$.\#
Then $m\bigl(\widetilde{\lm^{(H,x)}}\bigl)=\nu^{(H,x)}$.
\end{lemma}
\Dok
Let $Q$, $R$, $S_1$, $S_2$ be as in Definition~\ref{o:4:2} and
$\Nu=\Nu^{(H,x)}$, $\nu=\nu^{(H,x)}$.
By the main result of~\cite{Kleshchev_mullpaper} and
Lemma~\ref{l:4:2}, it suffices to prove that if $R>1$ then
\begin{equation}\label{eq:4:1.25}
G_p(\nu)=
\binom{2p}{2p-H}^Q
\binom{p+R-1}{p+R-H}
\binom{p}{p-H+1}^{x-2(Q+1)}
\binom{p-R+1}{p-H+1}
\end{equation}
and that if $R=1$ then
\begin{equation}\label{eq:4:1.375}
G_p(\nu)=
\binom{2p}{2p-H}^Q
\binom{p}{p-H+1}^{x-2Q}.
\end{equation}

Let us see first how the sets
$\<H,S_1,y\>$ and $\<p+H-R,S_2,x-2Q+y\>$
are situated with respect to one another.
Denote by $u_y$ the greatest element of the former set
and by $v_y$ the smallest element of the latter set.
We have $u_y=a^{(H,x)}_{y+p-2}$ and $v_y=b^{(H,x)}_{x-2Q+y}$.

For $y\ge0$, it follows from~(\ref{eq:4:1}) that
\begin{equation}\label{eq:4:1.5}
\begin{array}{l}
\ds
v_{y-1}-u_y=b^{(H,x)}_{x-2Q+y-1}-a^{(H,x)}_{y+p-2}=\\[12pt]
\ds
x-2Q+(H-R)\left[\tfrac{x-2Q+y-1}{p-H+1}\right]+
(R-1)\left[\tfrac{x-2Q+y-2}{p-H+1}\right]-
\left[\tfrac{R+y+p-3}{p-1}\right]\ge\\[16pt]
\ds
x-2Q-1+\left[\tfrac{x-2Q+y+p-H}{p-H+1}\right]
-\left[\tfrac{R+y+p-3}{p-1}\right]+
(R-1)\left[\tfrac{x-2Q+y-2}{p-H+1}\right].
\end{array}
\end{equation}

Since $x-2Q\ge R$, the first and the last summands of the last sum
are nonnegative for any $y\ge0$.
Moreover $p-H+1\le p-1$.
If $Q>0$ then we have $x-2Q+y+p-H=(H-2)Q+R+y+p-H>p+y+R-3$ and
therefore
\begin{equation}\label{eq:4:1.75}
\left[\tfrac{x-2Q+y+p-H}{p-H+1}\right]\ge\left[\tfrac{R+y+p-3}{p-1}\right].
\end{equation}

For $Q=y=0$, we get
$$
\left[\tfrac{x+p-H}{p-H+1}\right]\ge1\ge\left[\tfrac{R+p-3}{p-1}\right].
$$

Thus we have proved that $v_{y-1}\ge u_y$ for $0\le y\le Q$.
Recall that if $y\ge0$ then
being moved one position up
the greatest bead of $\<p+H-R,S_2,x-2Q+y\>$
takes position $v_{y-1}$.

Now one can clearly see that
$$
\vp^j(\Nu)=
\overline{(-\infty,H)\cup\<H,S_1,Q-j\>\cup\<p+H-R,S_2,x-Q-j\>}
\quad
\mbox{ for }
\quad
0\le j\le Q.
$$

Hence $h(\vp^j(\nu))=2p-H$ and $e( \vp^j(\nu) )=2p$ for $0\le j<Q$.
Moreover $h(\vp^Q(\nu))=p+R-H$.

If $R>1$ then $x-2Q\ge2$ and
$$
\vp^{Q+j}(\Nu)=\overline{(-\infty,p+H-1)\cup\<p+H,S_2,x-2Q-1-j\>}
\quad
\mbox{ for }
\quad
1\le j\le x-2Q-1.
$$
Hence $e(\vp^Q(\nu))=p+R-1$,
$e(\vp^{Q+j}(\nu))=p$ for $1\le j<x-2Q-1$,
$h(\vp^{Q+j}(\nu))=p-H+1$ for $1\le j\le x-2Q-1$,
$e(\vp^{x-Q-1}(\nu))=p-R+1$ and $\vp^{x-Q}(\nu)=\ems$.

If $R=1$ then
$$
\vp^{Q+j}(\Nu)=\overline{(-\infty,p+H-1)\cup\<p+H-1,S_2,x-2Q-j\>}
\quad
\mbox{ for }
\quad
0\le j\le x-2Q.
$$
Hence $e(\vp^{Q+j}(\nu))=p$ and $h(\vp^{Q+j}(\nu))=p-H+1$
for $0\le j<x-2Q$ and $\vp^{x-Q}(\nu)=\ems$.
\enddok

\begin{lemma}\label{l:4:4}
Let $h$, $i$, $x$ be integers such that
$0<i\le h<p$ and $x\ge i$.
If $i\ge h/2$ and $x\ge h$ then
$$
G_p(\lm^{(h,i,x)})=
\binom{A_0\cdots A_{h-i-1}}{R_0\cdots R_{h-i-1}}
\binom{p}{i}^{x-2(h-i)}
\binom
{A_0-2i\cdots A_{h-i-1}-2i}
{\eqn{R_0}{A_0}-\eqn{i}{2i}\cdots\eqn{R_{h-i-1}}{A_{h-i-1}}-\eqn{i}{2i}},
$$
if $i<h/2$ and $x\ge2i$ then
$$
G_p(\lm^{(h,i,x)})=
\binom{A_0\cdots A_{i-1}}{R_0\cdots R_{i-1}}
\binom{p}{i}^{x-2i}
\binom{A_0-2(h-i)\cdots A_{i-1}-2(h-i)}
{\eqn{R_0}{A_0}-\eqn{(h-i)}{2(h-i)}\cdots\eqn{R_{i-1}}{A_{i-1}}-\eqn{(h-i)}{2(h-i)}},
$$
if $x<2i,h$ then
$$
G_p(\lm^{(h,i,x)}){=}
\binom{A_0\cdots A_{x-i-1}}{R_0\cdots R_{x-i-1}}
\!
{\binom{p}{i+h-x}\!\!}^{2i-x}
\!
\binom{A_0-2(i+h-x)\cdots A_{x-i-1}-2(i+h-x)}
{\eqn{R_0}{A_0}-\eqn{(i+h-x)}{2(i+h-x)}\cdots\eqn{R_{x-i-1}}{A_{x-i-1}}-\eqn{(i+h-x)}{2(i+h-x)}}\!,
$$
where $A_j=p+h-1-2j$ and $R_j=h-j$,
\end{lemma}
\Dok
Let $S=[0,i-1]$, $\Lm=\Lm^{(h,i,x)}$ and $\lm=\lm^{(h,i,x)}$.

Case $i\ge h/2$ and $x\ge h$.
For $0\le j<h-i$, we have
$x-j-1\ge x-h+i\ge i$ and therefore the first element of
the set $\<0,S,x-j-1\>$ is
not less than $p$.
Hence
$$
\vp^j(\Lm)=\overline{(-\infty,j)\cup\<0,S,x-j\>\cup[i,h-j)}
\quad
\mbox{ for }
\quad
0\le j\le h-i.
$$
Therefore $h(\vp^j(\lm))=R_j$ for $0\le j\le h-i$ and
$e(\vp^j(\lm))=A_j$ for $0\le j<h-i$.
Next we have
$$
\vp^{h-i+j}(\Lm)=\overline{(-\infty,h-i)\cup\<h-i,S,x-2(h-i)-j\>}
\quad
\mbox{ for }
\quad
0\le j\le x-2(h-i).
$$
Hence
$e(\vp^{h-i+j}(\lm))=p$ for $0\le j<x-2(h-i)$,
$h(\vp^{h-i+j}(\lm))=i$ for $1\le j<x-2(h-i)$ and
$h(\vp^{x-(h-i)}(\lm))=h-i$.
Finally
$$
\vp^{x-(h-i)+j}(\Lm)=\overline{(-\infty,i+j)\cup[p,p+h-i-j)}
\quad
\mbox{ for }
\quad
0\le j\le h-i.
$$
Hence $e(\vp^{x-(h-i)+j}(\lm))=A_j-2i$ for $0\le j<h-i$ and
$h(\vp^{x-(h-i)+j}(\lm))=R_j-i$ for $1\le j\le h-i$.

Case $i<h/2$ and $x\ge2i$.
For $0\le j<i$, we have
$x-j-1\ge i$ and therefore the first element of the set
$\<0,S,x-j-1\>$ is
not less than $p$.
Hence we get
$$
\vp^j(\Lm)=\overline{(-\infty,j)\cup\<0,S,x-j\>\cup[i,h-j)}
\quad
\mbox{ for }
\quad
0\le j\le i.
$$
Therefore $h(\vp^j(\lm))=R_j$ and $e(\vp^j(\lm))=A_j$ for $0\le j<i$.
Next
$$
\vp^{i+j}(\Lm)=\overline{(-\infty,h-i)\cup\<p,S,x-2i-j\>}
\quad
\mbox{ for }
\quad
0\le j\le x-2i.
$$
Hence $e(\vp^{i+j}(\lm))=p$ and $h(\vp^{i+j}(\lm))=i$ for $0\le j<x-2i$.
Finally
$$
\vp^{x-i+j}(\Lm)=\overline{(-\infty,h-i+j)\cup[p,p+i-j)}
\quad
\mbox{ for }
\quad
0\le j\le i.
$$
Hence $e(\vp^{x-i+j}(\lm))=A_j-2(h-i)$ for $0\le j<i$ and
$h(\vp^{x-i+j}(\lm))=R_j-(h-i)$ for $0\le j\le i$.

Case $x<2i,h$.
For $ 0\le j<x-i$, we have $x-j-1\ge i$
and therefore the first element of
the set $\<0,S,x-j-1\>$ is
not less than $p$.
Hence
$$
\vp^j(\Lm)=\overline{(-\infty,j)\cup\<0,S,x-j\>\cup[i,h-j)}
\quad
\mbox{ for }
\quad
0\le j\le x-i
$$
Therefore $h(\vp^j(\lm))=R_j$ for $0\le j\le x-i$ and
$e(\vp^j(\lm))=A_j$ for $0\le j<x-i$.
Next
$$
\vp^{x-i+j}(\Lm)=\overline{(-\infty,x-i)\cup\<x-i,S,2i-x-j\>\cup[h,h-x+i)}
\quad
\mbox{ for }
\quad
0\le j\le 2i-x.
$$
Hence $e(\vp^{x-i+j}(\lm))=p$ for $0\le j<2i-x$ and
$h(\vp^{x-i+j}(\lm))=i+h-x$ for $1\le j<2i-x$.
Finally
$$
\vp^{i+j}(\Lm)=\overline{(-\infty,h-x+i+j)\cup[p,p+x-i-j)}
\quad
\mbox{ for }
\quad
0\le j\le x-i.
$$
Hence $e(\vp^{i+j}(\lm))=A_j-2(i+h-x)$ for $0\le j<x-i$ and
$h(\vp^{i+j}(\lm))=R_j-(i+h-x)$ for $0\le j\le x-i$.
\enddok

\begin{opr}\label{o:4:3}
Let $h$, $i$, $x$ be integers such that $0<i\le h<p$ and $x\ge i$.
Define
$$
\begin{array}{l}
\Mu^{(h,i,x)}=\overline{(-\infty,p)\cup[p+h-m,p+h)\cup
\<p,S,x\>},\\[6pt]
\mu^{(h,i,x)}=P(\Mu^{(h,i,x)}),
\end{array}
$$
where $m=\max\{i,i+h-x\}$ and $S=[0,h-m)\cup[h,p-1]$.
\end{opr}

Let $c^{(h,i,x)}_0$, $c^{(h,i,x)}_1$, \ldots
be all the elements of the set $\{n\in\Z:n\ge p, \rem(n,p)\in S\}$
written in ascending order.
An easy verification shows that
\begin{equation}\label{eq:4:2}
c^{(h,i,x)}_y=p+y+m+m\left[\tfrac{y+m-h}{p-m}\right].
\end{equation}

\begin{lemma}\label{l:4:5}
Let $h$, $i$, $x$ be integers such that
$0<i\le h<p$ and $x\ge i$.
Then $m(\lm^{(h,i,x)})=\mu^{(h,i,x)}$.
\end{lemma}
\Dok
Let $m$ and $S$ be as in Definition~\ref{o:4:3},
$A_j$ and $R_j$ as in Lemma~\ref{l:4:4} and
$\Mu=\Mu^{(h,i,x)}$, $\mu=\mu^{(h,i,x)}$, $S_j=p-j$.

Case $i\ge h/2$ and $x\ge h$.
Then $m=i$.
We have
$$
2p>A_0>\cdots>A_{h-i-1}>p>A_0-2i>\cdots>A_{h-i-1}-2i>0.
$$
Therefore by the main result of~\cite{Kleshchev_mullpaper}
and Lemma~\ref{l:4:4}, it suffices to prove that
$$
G_p(\mu)=
\binom{A_0\cdots A_{h-i-1}}{S_0\cdots S_{h-i-1}}
\binom{p}{p-i}^{x-2(h-i)}
\binom{A_0-2i\cdots A_{h-i-1}-2i}
{\eqn{S_0}{A_0}-\eqn{i}{2i}\cdots\eqn{S_{h-i-1}}{A_{h-i-1}}-\eqn{i}{2i}}.
$$

Note that $x-2(h-i)\ge0$.
Therefore for $0\le j<h-i$ we have $x-j-1\ge h-i$ and thus
the first element of the set $\<p,S,x-j-1\>$ is
not less than $p+h$.
Hence we have
$$
\vp^j(\Mu)=\overline{(-\infty,p+j)\cup[p+h-i,p+h-j)\cup\<p,S,x-j\>}
\quad
\mbox{ for }
\quad
0\le j\le h-i.
$$
Therefore $e(\vp^j(\mu))=A_j$ and $h(\vp^j(\mu))=S_j$ for $0\le j<h-i$.
Next
$$
\vp^{h-i+j}(\Mu)=\overline{(-\infty,p+i)\cup\<p+h,S,x-2(h-i)-j\>}
\quad
\mbox{ for }
\quad
0\le j\le x-2(h-i).
$$
Hence $e(\vp^{h-i+j}(\mu))=p$ and $h(\vp^{h-i+j}(\mu))=p-i$ for
$0\le j<x-2(h-i)$.
We have
$$
\vp^{x-(h-i)+j}(\Mu)=\overline{(-\infty,p+i+j)\cup[p+h,2p+h-i-j)}
\quad
\mbox{ for }
\quad
0\le j\le h-i.
$$
Hence
$e(\vp^{x-(h-i)+j}(\mu))=A_j-2i$ and
$h(\vp^{x-(h-i)+j}(\mu))=S_j-i$ for $0\le j<h-i$ and $\vp^x(\mu)=\ems$.

Case $i<h/2$ and $x\ge2i$.
We have
$$
2p>A_0>\cdots>A_{i-1}>p>A_0-2(h-i)>\cdots>A_{i-1}-2(h-i)>0.
$$
Therefore by the main result of~\cite{Kleshchev_mullpaper}
and Lemma~\ref{l:4:4}, it suffices to prove that
$$
G_p(\mu)=
\binom
{A_0\cdots A_{i-1}}
{S_0\cdots S_{i-1}}
\binom{p}{p-i}^{x-2i}
\binom
{A_0-2i\cdots A_{i-1}-2(h-i)}
{\eqn{S_0}{A_0}-\eqn{i}{2i}\cdots\eqn{S_{i-1}}{A_{i-1}}-\eqn{(h-i)}{2(h-i)}}.
$$

For $0\le j<i$, we have $x-j-1\ge h-m$.
Therefore the first element of the set $\<p,S,x-j-1\>$ is
not less than $p+h$.
Hence we get
$$
\vp^j(\Mu)=\overline{(-\infty,p+j)\cup[p+h-m,p+h-j)\cup\<p,S,x-j\>}
\quad
\mbox{ for }
\quad
0\le j\le i.
$$
Therefore $e(\vp^j(\mu))=A_j$ and $h(\vp^j(\mu))=S_j$ for $0\le j<i$.
For $0\le j\le x-2i$, we have
$$
\vp^{i+j}(\Mu){=}\left\{
\begin{array}{ll}
\!\overline{(-\infty,p+i)\cup[p+x-i-j,p+h-i)\cup[p+h,2p+x-i-j)}&\!\!\mbox{if }x<h;\!\!\\[6pt]
\!\overline{(-\infty,p+i)\cup\<p+i,S,x-2i-j\>}&\!\!\mbox{if }x\ge h.\!\!
\end{array}
\right.
$$
Hence $e(\vp^{i+j}(\mu))=p$ and $h(\vp^{i+j}(\mu))=p-i$ for $0\le j<x-2i$.
Finally
$$
\vp^{x-i+j}(\Mu)=\overline{(-\infty,p+h-i+j)\cup[p+h,2p+i-1-j]}
\quad
\mbox{ for }
\quad
0\le j\le i.
$$
Hence $e(\vp^{x-i+j}(\mu))=A_j-2(h-i)$, $h(\vp^{x-i+j}(\mu))=S_j-(h-i)$
for $0\le j<i$ and $\vp^x(\mu)=\ems$.

Case $x<2i,h$.
Then $m=i+h-x$.
We have
$$
2p>A_0>\cdots>A_{x-i-1}>p>A_0-2(i+h-x)>\cdots>A_{x-i-1}-2(i+h-x)>0.
$$
Therefore by the main result of~\cite{Kleshchev_mullpaper}
and Lemma~\ref{l:4:4}, it suffices to prove that
$$
G_p(\mu){=}
\ds\binom{A_0\cdots A_{x-i-1}}{S_0\cdots S_{x-i-1}}
\!\!
{\binom{p}{p-(i+h-x)}\!\!}^{2i-x}
\!\!\!
\binom
{A_0-2(i+h-x)\cdots A_{x-i-1}-2(i+h-x)}
{\eqn{S_0}{A_0}-\eqn{(i+h-x)}{2(i+h-x)}\cdots\eqn{S_{x-i-1}}{A_{x-i-1}}-\eqn{(i+h-x)}{2(i+h-x)}}\!.
$$
For $0\le j<x-i$, we have $x-j-1\ge i>x-i$ and therefore
the first element of the set $\<p,S,x-j-1\>$ is
not less than $p+h$.
Hence we get
$$
\vp^j(\Mu)=\overline{(-\infty,p+j)\cup[p+x-i,p+h-j)\cup\<p,S,x-j\>}
\quad
\mbox{ for }
\quad
0\le j\le x-i.
$$
Therefore $e(\vp^j(\mu))=A_j$ and $h(\vp^j(\mu))=S_j$ for $0\le j<x-i$.
Next
$$
\vp^{x-i+j}(\Mu)=\overline{(-\infty,p+i+h-x)\cup\<p+h,S,2i-x-j\>}
\quad
\mbox{ for }
\quad
0\le j\le 2i-x.
$$
Hence $e(\vp^{x-i+j}(\mu))=p$ and $h(\vp^{x-i+j}(\mu))=p-(i+h-x)$ for
$0\le j<2i-x$.
Finally
$$
\vp^{i+j}(\Mu)=\overline{(-\infty,p+i+h-x+j)\cup[p+h,2p+x-i-1-j]}
\quad
\mbox{ for }
\quad
0\le j\le x-i.
$$
Hence
$e(\vp^{i+j}(\mu))=A_j-2(i+h-x)$ and
$h(\vp^{i+j}(\mu))=S_j-(i+h-x)$
for $0\le j<x-i$ and $\vp^x(\mu)=\ems$.
\enddok

It is interesting to look at the partitions $\nu(2,x)$, when
$p>2$ and $x$ is an odd number greater than $2$.
By Lemmas~\ref{l:4:3} and~\ref{l:4:5}, we get
$m\bigl(\widetilde{\lm^{(2,x)}}\bigl)=\nu^{(2,x)}$ and
$m(\lm^{(2,x)})=\mu^{(2,2,x)}=\lm^{(p-2,x)}$.
We have $h(\nu^{(2,x)})=2p-2$ and $h(\lm^{(p-2,x)})=p-2$.
Hence
$\nu^{(2,x)}\not\ged\lm^{(p-2,x)}$ and
$(\nu^{(2,x)})^t\not\led(\lm^{(p-2,x)})^t$.
Since simple $\Sigma_{px}$-modules are self-dual,
by Theorem~\ref{t:cs:1} (or~\cite[Theorem~3.5(iv)]{Kleshchev_adacor}),
\cite[Theorem~4.4(b)]{Kleshchev_nak5}
and Propositions~\cite[II.2.14(4)]{Jantzen1},~\cite[2.1f]{Donkin3}, we have
$$
\begin{array}{l}
K\cong
\Ext_{\Sigma_{px}}\Bigl(D^{\widetilde{\lm^{(2,x)}}},D^{\lm^{(2,x)}}\Bigl)
\cong
\Ext_{S(N,px)}\Bigl(L\bigl((\nu^{(2,x)})^t\bigl),L\bigl((\lm^{(p-2,x)})^t\bigl)\Bigl)\\[12pt]
\cong
\Hom_{S(N,px)}\Bigl(\rad\Delta\bigl((\nu^{(2,x)})^t\bigl),L\bigl((\lm^{(p-2,x)})^t\bigl)\Bigl),
\end{array}
$$
where $N\ge px$ and $S(N,px)$ denotes the Schur algebra.
Definition~\ref{o:4:2} and~(\ref{eq:4:1}) show that
$\nu^{(2,x+2(p-1))}=\nu^{(2,x)}+(p^{2p-2})$.
Applying Corollary~5(4),
Theorem~4(a) and Lemma~4 from~\cite{Shchigolev9},
we obtain a negative solution of~Problem~2 from~\cite{Shchigolev9}
for the following values of the parameters:
$\lm:=\nu^{(2,3)}$, $q_i:=pi$, $n:=2p-2$, $V_i:=\rad S^{\lm+(q_i^n)}$.

In the remaining part of this section, we fix integers $H$, $x$, $i$
such that $1<H<p$, $H\nmid x$ and $0<i\le x,h$, where $h=H-\rem(x,H)$.
Let also $Q=\quo(x,H)$, $R=\rem(x,H)$,
$m=\max\{i,i+h-x\}$, $S_1=[0,p-1]\setminus\{h\}$,
$S_2=\{h\}\cup [H,p-1]$ and $S=[0,h-m)\cup[h,p-1]$.

\begin{lemma}\label{l:4:6}
${}$\\
{\rm (a)} if $c^{(h,i,x)}_y\ge b^{(H,x)}_z$ and $y,z>0$, then $c^{(h,i,x)}_{y-1}\ge b^{(H,x)}_{z-1}$;\\
{\rm (b)} if $a^{(h,i,x)}_y\le z$ and $y>0$, then $a^{(h,i,x)}_{y-1}\le z-1$;\\
{\rm (c)} if $a^{(h,i,x)}_y>z$, and $y\ge m-1$, $0\le m\le p$, then
         $a^{(h,i,x)}_{y-j}\ge z-j$ for $0\le j<m$.
\end{lemma}
{\bf Proof} follows from the mutual situation of the sets
$S$, $S_1$, $S_2$.
\enddok

\begin{teo}\label{t:4:1}
The following conditions are equivalent:
\begin{enumerate}
\item\label{t:4:1:p:1}
$m\bigl(\widetilde{\lm^{(H,x)}}\bigl)\ledeq m(\lm^{(h,i,x)})$;
\item\label{t:4:1:p:2}
$
\begin{array}{l}
\ds m\left[\tfrac{p-h+x-1}{p-m}\right]\ge H-Q-1+h\left[\tfrac{x-Q-1}{p-H+1}\right]+(R-1)\left[\tfrac{x-Q-2}{p-H+1}\right].
\end{array}
$
\end{enumerate}
\end{teo}
\Dok
Since $\shift(\Nu^{(H,x)})=\shift(\Mu^{(h,i,x)})=2p$
(we added one more row in the definition of $\Mu^{(h,i,x)}$
just to ensure this equality),
the first formula of~(\ref{eq:notation:1})
and Lemmas~\ref{l:4:3} and~\ref{l:4:5} show that condition~\ref{t:4:1:p:1}
of the current theorem is equivalent to $u\ledeq v$, where
$$
\begin{array}{l}
u=(b^{(H,x)}_{x-Q+p-H},\ldots,b^{(H,x)}_{x-Q},%
a^{(H,x)}_{Q+p-2},\ldots,a^{(H,x)}_Q,H-1\ldots,0),\\[12pt]
v=(c^{(h,i,x)}_{x+p-m-1},\ldots,c^{(h,i,x)}_x,p+h-1,\ldots,p+h-m,p-1,\ldots,0).
\end{array}
$$

It is easy to see that condition~\ref{t:4:1:p:2}
of the current theorem is equivalent to $u_1\le v_1$.
Therefore~\ref{t:4:1:p:1} implies~\ref{t:4:1:p:2}.

Now suppose that condition~\ref{t:4:1:p:2} is satisfied.
Note that $m\le h<H$ and thus $p-m\ge p-H+1$.
By the equivalence of the previous paragraph and Lemma~\ref{l:4:6}(a),
we get $u_j\le v_j$ for $1\le j\le p-H+1$.

Now let $p-H+1<j\le p-m$.
By~(\ref{eq:4:1}) and~(\ref{eq:4:2}), we have
$$
\begin{array}{l}
\ds v_j-u_j=c^{(h,i,x)}_{x+p-m-j}-a^{(H,x)}_{Q+2p-H-j}=
x-Q+m\left[\tfrac{x+p-j-h}{p-m}\right]-\left[\tfrac{2p-1+Q-j-h}{p-1}\right]\\[16pt]
\ds =x-Q-1+m\left[\tfrac{x+p-j-h}{p-m}\right]-\left[\tfrac{p+Q-j-h}{p-1}\right].
\end{array}
$$
Since $m\ge1$, $p-m\le p-1$ and $(x+p-j-h)-(p+Q-j-h)=x-Q>0$, we have
$$
m\left[\tfrac{x+p-j-h}{p-m}\right]\ge\left[\tfrac{p+Q-j-h}{p-1}\right]
$$
and $v_j\ge u_j$.

Case 1: $m=0$ or $m>0$, $v_{p-m+1}\ge u_{p-m+1}$.
By Lemma~\ref{l:4:6}(b), we get $v_j\ge u_j$ for $p-m<j\le p$.
Thus we have proved $v_j\ge u_j$ for $1\le j\le p$.
We have $v_j\le u_j$ for $p<j\le2p$,
since $v_j$ is an improper bead of $\Mu^{(h,i,x)}$ for such $j$.
Hence $\sigma_j(v)-\sigma_j(u)\ge \sigma_{2p}(v)-\sigma_{2p}(u)=0$
for any $p\le j\le2p$.
Therefore $u\ledeq v$.

Case 2: $m>0$ and $v_{p-m+1}<u_{p-m+1}$.
By Lemma \ref{l:4:6}(c), we get
$v_j\le u_j$ for $p-m<j\le p$.
We have $v_j\le u_j$ for $p<j\le2p$, since $v_j$ is an improper bead
of $\Mu^{(h,i,x)}$ for such $j$.
Hence $\sigma_j(v)-\sigma_j(u)\ge \sigma_{2p}(v)-\sigma_{2p}(u)=0$
for any $p-m\le j\le2p$.
Therefore $u\ledeq v$.
\enddok

{\bf Remark.}
The only fact we will need is that~\ref{t:4:1:p:1} implies~\ref{t:4:1:p:2}.
The reverse implication has been proved only to show
the impossibility of improving the bounds by replacing
a simpler condition~\ref{t:4:1:p:2}
with a more complicated condition~\ref{t:4:1:p:1}.

\section{Auxiliary upper bound}\label{auxest}

\subsection{Systems}\label{systems}
Introduce the following {\it staircase} abaci and partitions.
Let $k\ge0$ and $0<r_2<\cdots<r_k<p$, $i_1>\cdots>i_{k-1}\ge0$
be some integers.
We put
$$
\begin{array}{l}
\ds
\St(r_2,\ldots,r_k;i_1, \ldots,i_{k-1})=\\[10pt]
\ds
\overline{(-\infty,0)\cup
[pi_1,pi_1+r_2)\cup\
[pi_2+r_2,pi_2+r_3)\cup\cdots
[pi_{k-1}+r_{k-1},pi_{k-1}+r_k)},\\[10pt]
\ds
\st(r_2,\ldots,r_k;i_1, \ldots,i_{k-1})=
P(\St(r_2,\ldots,r_k;i_1, \ldots,i_{k-1})).
\end{array}
$$
For $k=1$, we assume $\St(\ems;\ems)=\overline{(-\infty,0)}$ and
$\st(\ems;\ems)=\ems$.

We have already met special cases of staircase abaci in
\sectsign~\ref{mullmap}.
Indeed, let $h$, $i$, $x$ be integers such that $0<i\le h<p$ and $x\ge0$.
We put $q=\quo(x,i)$ and $r=\rem(x,i)$.
Then
\begin{equation}\label{eq:auxest:0.5}
\Lm^{(h,i,x)}=\left\{
\begin{array}{ll}
\St(i,h;q,0)&\mbox{ if }r=0;\\
\St(r,i,h;q+1,q,0)&\mbox{ if }r>0.
\end{array}
\right.
\end{equation}

\begin{lemma}\label{l:auxest:1}
Every solution of the system
\begin{equation}\label{eq:auxest:1}
\left\{
\begin{array}{l}
h(\lm)<p,\\
\core(\lm)=\ems,\\
\res A=\res B
\text{ for any $\lm$-normal nodes $A$ and $B$}.
\end{array}
\right.
\end{equation}
has the form
$\st(r_2,\ldots,r_k;i_1, \ldots,i_{k-1})$.
\end{lemma}
\Dok
Clearly, these partitions satisfy the system.
Prove the converse fact by induction on the number $n$ of nodes in $[\lm]$.
This is obviously true for $n=0$.
Now let $n>0$ and $\lm$ be a solution of system~(\ref{eq:auxest:1}).
Since $\lm\ne\ems=\core(\lm)$, there is at least one
rim $p$-hook of $\lm$.
Denote by $\bar\lm$ the partition obtained from $\lm$
by removing the highest of these hooks.
By Lemma~\ref{l:2:1}, the partition $\bar\lm$
also satisfies system~(\ref{eq:auxest:1}).
By the inductive hypothesis,
$\bar\lm=\st(r_2,\ldots,r_k;i_1, \ldots,i_{k-1})$.

Take the abacus $\Lm$ such that
$$\lm=P(\Lm)\mbox{ and }
\shift(\Lm)=\shift(\St(r_2,\ldots,r_k;i_1, \ldots,i_{k-1})).
$$
We see that there is a movable down bead $a$ of
$\St(r_2,\ldots,r_k;i_1, \ldots,i_{k-1})$ such that
$\Lm$ is obtained from $\St(r_2,\ldots,r_k;i_1, \ldots,i_{k-1})$
by moving $a$ one position down.

Case $k=1$.
Since $h(\lm)<p$, we have $-p<a<0$ and
$\lm=\st(1,-a;1,0)$ for $a<-1$ and $\lm=\st(1;1)$ for $a=-1$.

Case $k>1$.
We can assume $i_1>0$,
as otherwise $\bar\lm=\ems$ and we are
under the conditions of the previous case.
Since $h(\lm)<p$, we have $r_k-p<a$.
If $a<-1$ or $a\le-1$ and $i_k>1$,
then $\Lm$ contains the following two normal beads
in different runners: $pi_1$ and $a+p$.
Therefore, the third equation of system~(\ref{eq:auxest:1})
is violated for $\lm$ and this case is impossible.
If $a=-1$ and $i_k=1$, then $\lm=\st(r_2+1,\ldots,r_k+1;i_1, \ldots,i_{k-1})$
($k=2,3$).

Now suppose that $a\ge0$.
Since the third equation of system~(\ref{eq:auxest:1})
holds for $\lm$, $a$ can take only the following values:
$pi_t+r_t$, where $t=1$, \ldots, $k-1$ and we assume $r_1=0$.
Otherwise $pi_1$ and $a+p$ would be normal beads of $\Lm$
from different runners.
Now directly from the definition of staircase abacus
one can see that $\lm$ has the desired form.
\enddok

Let us calculate $e(\st(r_2,\ldots,r_k;i_1, \ldots,i_{k-1}))$.
It is zero for $k=1$.
Therefore consider the case $k>1$.
Define the sequence
$1=a_1<\cdots<a_l\le k-1$ by the following rule:
$a_{j+1}=a_j+1$ if $a_j+1\le k-1$ and $i_{a_j+1}<i_{a_j}-1$;
$a_{j+1}=a_j+2$ if $a_j+2\le k-1$ and $i_{a_j+1}=i_{a_j}-1$.
From the definition of \sectsign~\ref{abaci}, we get
$$
e(\st(r_2,\ldots,r_k;i_1, \ldots,i_{k-1}))=
\left\{
\begin{array}{ll}
pl& \mbox{ if }i_{a_l}>0;\\[6pt]
p(l-1)+r_k-1&\mbox{ if }i_{a_l}=0.
\end{array}
\right.
$$
Therefore if $e(\lm)<2p$ and $\lm$ is a staircase partition,
then $\lm$ has one of the following forms:
$\ems$, $\st(r_2,r_3;i_1,0)$, $\st(r_2,r_3,r_4,;i_1,i_1-1,0)$.
This fact and~(\ref{eq:auxest:0.5}) yield

\begin{lemma}\label{l:auxest:2}
Every solution of the system
\begin{equation}\label{eq:auxest:2}
\left\{
\begin{array}{l}
h(\lm)<p,\\
\core(\lm)=\ems,\\
\res A=\res B
\text{ for any $\lm$-normal nodes $A$ and $B$},\\
h(\lm)+h(m(\lm))<2p.
\end{array}
\right.
\end{equation}
has the form $\lm^{(h,i,x)}$,
where $0<i\le h<p$ and $x\ge0$.
\end{lemma}
The last inequalities can be strengthen by $i\le x$.

\subsection{Bound}\label{est}
We state the following known result.

\begin{utv}\label{u:est:1}
Let $n>2$ and $V$ be a $K\Sigma_{n-1}$-module.
Then
\begin{enumerate}
\itemsep=0pt
\renewcommand{\labelenumi}{{\rm \theenumi}}
\renewcommand{\theenumi}{(\alph{enumi})}
\item\label{u:est:1:c:1} $V\uparrow^{\Sigma_n}\otimes\sgn_n\cong(V\otimes\sgn_{n-1})\uparrow^{\Sigma_n}$;
\item\label{u:est:1:c:2} $(\Ind^\alpha V)\otimes\sgn_n\cong\Ind^{-\alpha}(V\otimes\sgn_{n-1})$.
\end{enumerate}
\end{utv}
\Dok\ref{u:est:1:c:1}
The isomorphism is given by
$(\sigma_i\otimes v)\otimes u\mapsto\sigma_i\otimes (v\otimes u)$,
where $\sigma_1$, \ldots, $\sigma_n$ are representatives of
the left cosets of $\Sigma_n$ over $\Sigma_{n-1}$ having the same sign,
$v\in V$ and $u$ is a basis of the sign representation of $\Sigma_n$.

\ref{u:est:1:c:2}
follows from \ref{u:est:1:c:1} and
the arguments from the proof of~\cite[Theorem~4.7]{Kleshchev3}
applied to the induction operator instead of
the restriction operator.
\enddok

\begin{opr}\label{o:auxest:1}
Let $\pi(H,x,i)$ be satisfied if and only if
\begin{enumerate}
\item\label{pidef:1} $H,x,i\in\Z$, $2<H<p$,
$H\nmid x$, $0<i\le x,h$ and $x>2$;
\item\label{pidef:2}
$
\ds
m\left[\tfrac{p-h+x-1}{p-m}\right]{\ge}H{-}Q{-}1+h\left[\tfrac{x-Q-1}{p-H+1}\right]{+}(R{-}1)\left[\tfrac{x-Q-2}{p-H+1}\right],
$
\end{enumerate}
where $Q=\quo(x,H)$, $R=\rem(x,H)$, $h=H-R$ and $m=\max\{i,i+h-x\}$.
\end{opr}

\begin{opr}\label{o:auxest:2}
Let $H,x,i$ be integers for which $\pi(H,x,i)$ is satisfied.
Denote by $\epsilon(H,x,i)$ the sequence such that
$-R+\Lm^{(h,i,x)}=\H_{\epsilon(H,x,i)}(\Lm^{(H,x)})$,
where $R=\rem(x,H)$ and $h=H-R$.
\end{opr}

Clearly the required sequence $\epsilon(H,x,i)$ exists
and is given by
$$
\epsilon(H,x,i)=\bigl(
(-Q-1)^R,(-Q)^{h-i},(q-Q)^{i-r},(q+1-Q)^r
\bigl),
$$
where $Q=\quo(x,H)$, $R=\rem(x,H)$, $h=H-R$,
    $q=\quo(x,i)$, $r=\rem(x,i)$.

\begin{opr}\label{o:auxest:3}
A pair $(\nu,\mu)$ is called minimal if
\begin{itemize}
\itemsep=0pt
\item $\nu$ is an almost completely splittable partition;
\item $\mu$ is a $p$-regular partition;
\item $\nu\not\ged\mu$ and $\nu\sim\mu$;
\item for any $\mu$-good node $A$ there exists a $\nu$-good node $B$
      such that $\res B=\res A$ and $\nu_B$ is not an
      almost completely splittable partition.
\end{itemize}
\end{opr}

It follows directly from the definition that
if $B$ is a $\nu$-normal node,
$\nu$ is an almost completely splittable partition and $\nu_B$ is not,
then $\chi(\lm)=p$ and $B=(1,\nu_1)$,
where $\lm$ is the preimage of $\nu$.
Therefore if $(\nu,\mu)$ is a minimal pair,
then all $\mu$-normal nodes have residue
$\overline{\nu_1-1}=\overline{\lm_1}$.

\begin{lemma}\label{l:auxest:3}
Let $(\nu,\mu)$ be a minimal pair of partitions of $n$.
We put $x=n/p$ and $H=h(\lm)$, where $\lm$ is the preimage of $\nu$.
Suppose that $(\nu,\mu)$ is different from $((p^2,p^2-p),(2p^2-p))$
and from $(\widetilde{\lm^{(H,2)}},\lm^{(H-2,2)})$
and that $\Ext_{\Sigma_n}(D^\nu,D^\mu)\ne0$.
Then $\nu=\widetilde{\lm^{(H,x)}}$ and
$\mu=\lm^{(H-\rem(x,H),i,x)}$, where $\pi(H,x,i)$ is satisfied.
\end{lemma}
\Dok
Let $\alpha$ be the residue of all $\mu$-normal nodes and
$A$ be a $\nu$-good node of residue $\alpha$.
We have $n\ge p>2$.
Since $h(\nu)<p$, by Proposition~\ref{u:notation:0.5}, we have
$$
\Hom_{\Sigma_n}(\rad S^\nu,D^\mu)\cong\Ext_{\Sigma_n}(D^\nu,D^\mu)\ne0.
$$
Hence $\nu\led\mu$ and in particular $h(\mu)\le h(\nu)$.
Let $\nu=\tilde\lm$, where $\lm$ is a big partition.
Since $\nu_A$ is not almost completely splittable, we get
$\chi(\lm)=p$ and $A$ is in the first row.

We have $h(\nu)\le h(\lm)$.
Suppose $h(\nu)<h(\lm)$.
Then $\nu$ is completely splittable.
By the hypothesis of the current theorem and Theorem~\ref{t:cs:1},
we get that $\nu$ is a big partition and $\mu=\tilde\nu$.
The rightmost node $A'$ of the first row of $[\mu]$ is removable and
therefore is normal.
However $\res A'$ is equal to the residue of
the rightmost node of the last row of $[\nu]$,
which is distinct from $\res A$, as
$\nu$ is completely splittable and $h(\nu)>1$.
The resulting contradiction gives $h(\nu)=h(\lm)$.

By Corollary~\ref{sl:2:0}, the residue of any $\core(\mu)$-normal node,
and by $\lm\sim\mu$ also of any $\core(\lm)$-normal node, is $\res A$,
which in turn equals the residue of the bottom $\lm$-removable node.
By Lemma~\ref{l:3:1}, we have $\lm=\lm^{(H,x)}$, $1<H<p$, $x>0$ and $H\nmid x$.
The case $x=1$ is impossible,
as we would have a contradiction $h(\nu)<h(\lm)$.
Thus $x>1$.
By Lemma~\ref{l:3:2}, we get $h(\mu)<h(\lm)=H$.
We put $Q=\quo(x,H)$, $R=\rem(x,H)$ and $h=H-R$.

If $R>1$ we put $b=(Q+1)p$ and if $R=1$ we put $b=Qp$.
Since $x>1$, we get that $b$ is the only normal
and thus good bead of $\widetilde{\Lm^{(H,x)}}$ belonging to runner zero .
Denote by $B$ the node of $\nu$ corresponding to $b$.
By the second formula of~(\ref{eq:notation:1}), we have
$\res B=-\overline{\shift(\Lm^{(H,x)})}=-\bar H$.
By Lemma~\ref{l:inequal:3}, we have
$$
\begin{array}{l}
\ds0<\dim\Ext_{\Sigma_n}(D^\nu,D^\mu)\le\\[10pt]
\ds\dim\Ext_{\Sigma_{n-1}}(D^{\nu_B},\Res_{-\bar H}D^\mu)+
\dim\Hom_{\Sigma_{n-1}}(\rad\Ind^{-\bar H}D^{\nu_B},D^\mu).
\end{array}
$$
Since $\alpha=\bar R-\bar H\ne-\bar H$, the hypothesis of the
current theorem yields $\Res_{-\bar H}D^\mu=0$.
Hence $\Hom_{\Sigma_{n-1}}(\rad\Ind^{-\bar H}D^{\nu_B},D^\mu)\ne0$ and
in particular $[\Ind^{-\bar H}D^{\nu_B}:D^\mu]>0$.
Multiplying the modules from the last inequality by $\sgn_n$,
we get by Proposition~\ref{u:est:1}\ref{u:est:1:c:2} that
\begin{equation}\label{eq:est:2.5}
[\Ind^{\bar H}D^{m(\nu_B)}:D^{m(\mu)}]>0.
\end{equation}
It follows from the proof of~\cite[Theorem~4.7]{Kleshchev3} that
$m(\nu_B)=m(\nu)_{B^{(m)}}$, where $B^{(m)}$ is the $m(\nu)$-good node
of residue $\bar H$.
By Lemma~\ref{l:4:3}, we get $m(\nu)=P(\Nu^{(H,x)})$.
For brevity until the end of the current proof we shall
use the notation $\Nu=\Nu^{(H,x)}$, $a_i=a^{(H,x)}_i$, $b_i=b^{(H,x)}_i$
(see Definition~\ref{o:4:2} and~(\ref{eq:4:1})).
We shall prove that $m(\nu)\ledeq m(\mu)$.

Case 1: $R>1$.
We have $a_{Q+p-1}<b_{x-Q}$.
Indeed if $Q>0$, then
taking into account~(\ref{eq:4:1.75}) the substitution
$y=Q+1$ into~(\ref{eq:4:1.5}) yields $b_{x-Q}-a_{Q+p-1}\ge x-2Q-1\ge R-1>0$.
On the other hand, if $Q=0$ then $a_{Q+p-1}=H+p<b_x=b_{x-Q}$, since $x>1$.

Denote by $b^{(m)}$ the element of $\<p+H-R,S_2,x-Q\>$
from runner $H$.
Clearly, $b^{(m)}$ is a normal bead of $\Nu$.
The only initial bead of $\Nu$ distinct from $b^{(m)}$ and
belonging to runner $H$ is $a_Q$ in the case where $Q>0$ and $p-1\mid Q$.
However this bead is not normal, as in this case
$a_Q+p-1=a_{Q+p-2}$ and $a_Q+p=a_{Q+p-1}<b_{x-Q}$,
whence $N(a_Q+p-1)=1$ and $N(a_Q+p)=0$.
Therefore $b^{(m)}$ is a good bead of $\Nu$,
the node $B^{(m)}$ corresponds to this bead,
and $P(\Nu_{b^{(m)}})=m(\nu)_{B^{(m)}}=m(\nu_B)$.

All the initial spaces of $\Nu_{b^{(m)}}$ from runner $H$ are:
${b^{(m)}}$; $H$ if $Q>0$; $a_{Q+p-1}$ if $p-1\mid Q$.
If in the last case $Q>0$, then the space $a_{Q+p-1}$ is not conormal.
We put $c=H$ if $Q>0$ and $c=H+p$ if $Q=0$.
We have $c\le a_{Q+p-1}<b_{x-Q}\le b^{(m)}$.
Therefore all the conormal spaces of $\Nu_{b^{(m)}}$
from runner $H$ are ${b^{(m)}}$ and $c$.

Case 1.1: $m(\mu)\ne P(\Nu_{b^{(m)}}^c)$.
Denote by $C$ the addable node of $m(\nu_B)$
corresponding to the space $c$ of $\Nu_{b^{(m)}}$.
If we suppose $[S^{m(\nu)}:D^{m(\mu)}]=0$, then by~(\ref{eq:est:2.5}),
Lemma~\ref{l:spechtex:2.5} applied to the following parameters:
$$
\hbox to \textwidth
{$\ds r:=n-1$,\hfil $\ds \lm:=m(\nu_B)$,
\hfil $\ds \alpha:=\bar H$,\hfil $\ds k:=2$,\hfil $\ds B_1:=C$,\hfil $\ds B_2:=B^{(m)}$,\hfil
$\ds \gamma:=m(\mu)$}
$$
yields a contradiction $m(\mu)=m(\nu_B)^C=P(\Nu_{b^{(m)}}^c)$.
Thus $[S^{m(\nu)}:D^{m(\mu)}]>0$, $m(\nu)\ledeq m(\mu)$ and
$m(\nu)\ledeq m(\mu)$.

Case 1.2: $m(\mu)=P(\Nu_{b^{(m)}}^c)$.
If one of the conditions $R=2$ or $Q=0$ is violated,
then $\Nu_{b^{(m)}}^c$ contains the following normal beads not belonging
to runner $H-R$: $b^{(m)}-1$ if $R>2$;
the bead of $\<H,S_1,Q\>$ from runner $H-1$ if $R=2$ and $Q>0$.
Therefore the partition $m(\mu)=P(\Nu_{b^{(m)}}^c)$ contains
normal nodes of residue different from $\bar H-\bar R=-\alpha$,
which is a contradiction.

We have $x=2$, $h(\nu)=H>2$, $\nu=\widetilde{\lm^{(H,2)}}$ and
\begin{equation}\label{eq:est:2.75}
\!
\!
\!
\begin{array}{lcl}
\ds \Nu&=&\overline{(-\infty,H-2+p)\cup\{H-1+p\}\cup(H+p,2p)\cup\{H-2+2p,H+2p\}},\\[6pt]
\ds \Nu_{b^{(m)}}^c&=&\overline{(-\infty,H-2+p)\cup[H+p,2p)\cup\{H-2+2p,H-1+2p\}}.
\end{array}
\end{equation}
Hence $m(\mu)=(H^2,2^{p-H})=\lm^{(p-H+2,2,2)}$.
By Lemma~\ref{l:4:5}, we get
$\mu=\lm^{(H-2,2)}$, which contradicts the hypothesis of the theorem.

Case 2: $R=1$.
Since $x>1$, we have $Q>0$.
In the proof of Lemma~\ref{l:4:3}, it was shown that
$a_{Q+p-2}\le b_{x-Q-1}<b_{x-Q}$.
Denote by $b^{(m)}$ the element of $\<H,S_1,Q\>$ from runner $H$.
Since $Q>0$, this bead is normal in $\Nu$.
Clearly, there is no other normal bead of $\Nu$ in runner $H$.
Therefore $b^{(m)}$ is a good bead of $\Nu$,
the node $B^{(m)}$ corresponds to this bead,
and $P(\Nu_{b^{(m)}})=m(\nu)_{B^{(m)}}=m(\nu_B)$.

All the initial spaces of $\Nu_{b^{(m)}}$ from runner $H$ are:
${b^{(m)}}$; $H$; $b_{x-Q+p-H+1}$ if $p-H+1\mid x-Q-1$.
In the last case the space $b_{x-Q+p-H+1}$ is not conormal in $\Nu_{b^{(m)}}$,
since $b^{(m)}\le a_{Q+p-2}<b_{x-Q}$.
We put $c=H$.
Since $Q>0$, we have $c<b^{(m)}$.
Therefore all the conormal spaces of $\Nu_{b^{(m)}}$
from runner $H$ are ${b^{(m)}}$ and $c$.
Denote by $C$ the addable node of $m(\nu_B)$ corresponding
to the space $c$ of $\Nu_{b^{(m)}}$.

Case 2.1: $m(\mu)\ne P(\Nu_{b^{(m)}}^c)$
is similar to case~1.1.

Case 2.2: $m(\mu)=P(\Nu_{b^{(m)}}^c)$.
We have $h(m(\mu))=2p-H+1$ and $h(\mu)\le H-1$.
Hence $e(\Nu_{b^{(m)}}^c)\le h(\mu)+h(m(\mu))\le 2p$.
Clearly, this can happen only if $b^{(m)}=a_{Q+p-2}=H+p$
(arguments similar to the ones used in the proof of Lemma~\ref{l:4:3}).
In this case
\begin{equation}\label{eq:est:3}
\Nu_{b^{(m)}}^c=\overline{(-\infty,H-1)\cup[H,H+p)\cup\<p+H-1,S_2,x-Q\>}
\end{equation}
and $m(\mu)=P(\Nu_{b^{(m)}}^c)$ is a $p$-singular partition,
which is a contradiction.

Thus we have considered all possible cases and proved that
$m(\nu)\ledeq m(\mu)$.
We have $h(\mu)<h(\nu)<p$ and $h(m(\mu))\le h(m(\nu))$.
By~(\ref{eq:4:1.25}) and~(\ref{eq:4:1.375}), we have
$$
h(\mu)+h(m(\mu))<h(\nu)+h(m(\nu))=e(\nu)+[p\nmid e(\nu)]\le 2p.
$$
Therefore $\mu$ satisfies system~(\ref{eq:auxest:2})
and Lemma~\ref{l:auxest:2} yields $\mu=\lm^{(h',i,x)}$,
where $0<i\le h',x$ and $h'<p$.
Since the only normal bead of $\lm^{(h',i,x)}$
has residue $-\bar{h'}$, the hypothesis of the current theorem yields
$-\bar{h'}=\alpha=\bar R-\bar H$.
Since $0<H-R<p$, we have $-p<H-R-h'<p$.
But $H-R-h'$ is divisible by $p$, whence $h'=H-R=h$.

From $m(\nu)\ledeq m(\mu)$ and Theorem~\ref{t:4:1}, it follows that
part~\ref{pidef:2} of the condition $\pi(H,x,i)$ is satisfied.
We claim that part~\ref{pidef:1} of the condition $\pi(H,x,i)$
is satisfied.
It suffices to prove that $H,x>2$.

Indeed, if $H=2$ then
$h=1$, $x$ is an odd number greater than $2$ and
$\nu=(\tfrac{x+1}2p,\tfrac{x-1}2p)$, $\mu=(px)$.
By~\cite[Lemma~3.5(iv)]{Kleshchev_adacor},
we get $x=2p-1$ and $(\nu,\mu)=((p^2,p^2-p),(2p^2-p))$.
This is a contradiction.

If $x=2$ then part~\ref{pidef:2} of the condition $\pi(H,x,i)$
is not satisfied, contrary to what was proved.
\enddok

\begin{lemma}\label{l:auxest:4}
${}$\\[-12pt]
\begin{enumerate}
\itemsep=0pt
\renewcommand{\labelenumi}{{\rm \theenumi}}
\renewcommand{\theenumi}{(\alph{enumi})}
\item\label{l:auxest:4:p:1} $\Ext_{\Sigma_{2p^2-p}}(D^{(p^2,p^2-p)},D^{(2p^2-p)})\cong K$;
\item\label{l:auxest:4:p:2} $\Ext_{\Sigma_{2p}}(D^{\widetilde{\lm^{(H,2)}}},D^{\lm^{(H-2,2)}})\cong K$, where $2<H<p$.
\end{enumerate}
\end{lemma}
\Dok
\ref{l:auxest:4:p:1} follows from~\cite[Lemma~3.5(iv)]{Kleshchev_adacor}.

\ref{l:auxest:4:p:2}
It follows from~(\ref{eq:est:2.75}) and Lemma~\ref{l:4:3} that
$m\bigl(\widetilde{\lm^{(H,2)}}\bigl)=\widetilde{\lm^{(p-H+2,2)}}$.
Hence by Theorem~\ref{t:cs:1}, we get
\begin{equation}\label{eq:auxest:4}
K\cong\Ext_{\Sigma_{2p}}(D^{\lm^{(p-H+2,2)}},D^{m\bigl(\widetilde{\lm^{(H,2)}}\bigl)})
\cong\Ext_{\Sigma_{2p}}(D^{\widetilde{\lm^{(H,2)}}},D^{m(\lm^{(p-H+2,2)})}).
\end{equation}
By Lemma~\ref{l:4:5}, we have
$m(\lm^{(p-H+2,2)})=m(\lm^{(p-H+2,2,2)})=\lm^{(H-2,2)}$
(see case~1.2 of the proof of the previous lemma).
Substituting this value into~(\ref{eq:auxest:4}),
we get the desired equivalence.
\enddok

\section{Upper bound}\label{upperest}

\subsection{General construction}
To make the approach of \sectsign~\ref{cs_gen_const} applicable
to almost completely splittable partitions, we shall modify it as follows.

Let $X$ be a set satisfying the same conditions as in
\sectsign~\ref{cs_gen_const} and we somehow know
a map $\zeta:X\to\Z'$ such that
\begin{equation}\label{eq:upperest:1}
\dim\Ext_{\Sigma_n}(D^\nu,D^\mu)\le\zeta(\nu,\mu)
\mbox{ for any pair }(\nu,\mu)\in X\mbox{ of partitions of }n.
\end{equation}

Define the map $U:X\to\Z'$ by induction as follows.
We put $U(\ems,\ems)=0$.
Now let $(\nu,\mu)$ be a pair of nonempty partitions of $X$.
For any $\mu$-good node $A$, let $m_A(\nu,\mu)$ equal:
\begin{itemize}
\item $\varepsilon(\nu,\mu_A)$ if there is no $\nu$-good node
      of residue $\res A$;
\item $U(\nu_B,\mu_A)+\varepsilon(\nu,\mu_A)$ if
      there is a $\nu$-good node $B$ of residue $\res A$ and
      $(\nu_B,\mu_A)\in X$;
\item $+\infty$ if there is a $\nu$--good node $B$ of residue $\res A$ and
      $(\nu_B,\mu_A)\notin X$.
\end{itemize}
We put
$U(\nu,\mu)=\min\bigl(\{\zeta(\nu,\mu)\}\cup
\{m_A(\nu,\mu): A$ is a $\mu$-good node$\}\bigl)$.

\begin{lemma}\label{l:upperest:1}
Let $(\nu,\mu)\in X$.
Then $\dim\Ext_{\Sigma_n}(D^\nu,D^\mu)\le U(\nu,\mu)$,
where $\nu,\mu\vdash n$.
\end{lemma}
{\bf Proof} is by induction on $n$ applying Lemma~\ref{l:inequal:1}.
\enddok

\subsection{Case of almost completely splittable partitions}

Let us keep the following notation until the end of this section:
$$
X=\{(\nu,\mu):\nu\mbox{ is almost completely splittable},\;%
              \mu\mbox{ is $p$-regular},\;%
              \nu\not\ged\mu\;,%
              \nu\sim\mu\}.%
$$

Define $\zeta$ as follows.
Let $(\nu,\mu)\in X$ be a pair of partitions of $n$.
We put $x=n/p$ and $H=h(\lm)$, where $\lm$ is the preimage of $\nu$.
We define $\zeta(\nu,\mu)=+\infty$ except the following cases:
\begin{enumerate}
\itemsep=0pt
\item $\mu\not\ged\nu$. We put $\zeta(\nu,\mu)=0$.
\item $\mu\ged\nu$ and $(\nu,\mu)$ is minimal.
{
\renewcommand{\labelenumii}{\theenumi}
\renewcommand{\theenumii}{}
\renewcommand{\theenumi}{(\arabic{enumi}.\arabic{enumii})}
\begin{enumerate}
\itemsep=0pt
\item\label{zeta:c:2:1} $(\nu,\mu)=((p^2,p^2-p),(2p^2-p))$.
      We put $\zeta(\nu,\mu)=1$.
\item\label{zeta:c:2:2} $(\nu,\mu)=(\widetilde{\lm^{(H,2)}},\lm^{(H-2,2)})$,
      where $H>2$.
      We put $\zeta(\nu,\mu)=1$.
\item\label{zeta:c:2:3} Cases~\ref{zeta:c:2:2} and~\ref{zeta:c:2:3}
      do not hold and there is no $i$ such that
      $\nu=\widetilde{\lm^{(H,x)}}$, $\mu=\lm^{(H-\rem(x,H),i,x)}$ and
      $\pi(H,x,i)$ is satisfied.
      We put $\zeta(\nu,\mu)=0$.
\end{enumerate}}
\end{enumerate}
It follows from Proposition~\ref{u:notation:0.5} and
Lemmas~\ref{l:auxest:3},~\ref{l:auxest:4} that
property~(\ref{eq:upperest:1}) holds for $\zeta$ we have just defined.

\begin{teo}\label{t:upperest:1}
Let $\nu,\mu$ be $p$-regular partitions of $n$ such that
$\nu$ is almost completely splittable and $\nu\not\ged\mu$.
Then $\Ext_{\Sigma_n}(D^\nu,D^\mu)=0$ except the case when
$\mu=\H_\epsilon(\lm)$, where $\lm$ is the preimage of $\nu$, $H=h(\lm)$
and one of the following conditions holds:
\begin{enumerate}
\item\label{t:upperest:1:c:1} $H=2$ and $\epsilon=(1,-1),(-p,p)$;
\item\label{t:upperest:1:c:2} $H=3$ and $\epsilon=(0,-1,1),(-1,1,0),(-1,-1,2)$;
\item\label{t:upperest:1:c:3} $H>3$ and $\epsilon=(0,-1,0^{H-3},1),(-1,0^{H-3},1,0),(-1,-1,0^{H-4},1,1)$;
\item\label{t:upperest:1:c:4} $\epsilon=\epsilon(H,x,i)$ for $x$ and $i$
     such that $\pi(H,x,i)$ is satisfied.
\end{enumerate}
In cases~\ref{t:upperest:1:c:1}-\ref{t:upperest:1:c:3} the bound
$\Ext_{\Sigma_n}(D^\nu,D^\mu)\le1$ holds.
\end{teo}
\Dok
In view of Proposition~\ref{u:notation:0.75},
it suffices to assume $\nu\sim\mu$.
Therefore we must prove that $U(\nu,\mu)=0$
for any pair $(\nu,\mu)\in X$ except
cases~\ref{t:upperest:1:c:1}-\ref{t:upperest:1:c:4} and that
$U(\nu,\mu)\le1$ in cases~\ref{t:upperest:1:c:1}-\ref{t:upperest:1:c:3}.
Apply induction on $n$.
If $n=0$ then by definition we have $U(\nu,\mu)=0$.

Now let $n>0$ and suppose that theorem is true for partitions of numbers
less than $n$.
Choose some abaci $\Mu$ and $\Lm$ of the same shift
such that $\mu=P(\Mu)$ and $\lm=P(\Lm)$.

Case 1: $\mu\not\ged\nu$. Then $U(\nu,\mu)=\zeta(\nu,\mu)=0$.

Case 2: $\mu\ged\nu$ and the pair $(\nu,\mu)$ is not minimal.
Let $A$ be a $\mu$-good node such that either there is no
$\nu$-good node of residue $\res A$ or such a node $B$
exists and $\nu_B$ is almost completely splittable.
Denote by $a$ the bead of $\Mu$ corresponding to $A$.

First consider the case $\epsilon(\nu,\mu_A)=1$.
Then $\tilde\Lm=(\Mu_a)^c$ for some conormal space $c$ of $\Mu_a$.
Since $\nu\sim\mu$ and $\nu\led\mu$,
we get that $a$ and $c$ belong to the same runner
and moreover $a$ is below $c$.
We have $\Mu=(\tilde\Lm_c\vphantom{\Lm})^a$.

If $H=2$ then we get $c=\b^\Lm-p$, $a=\b^\Lm+p$ and $\b^\Lm=\b_\Lm+1$.
Hence $\Mu=\{\b_\Lm-p,\b^\Lm+p\}=\H_{(1,-1)}(\Lm)$.

If $H>3$ then either $c=\b^\Lm-p$, $a=\b^\Lm$, $\b^\Lm=\b^\Lm(2)+1$
or $c=\b_\Lm(2)$, $a=\b^\Lm(2)+p$, $\b_\Lm(2)=\b_\Lm+1$.
In the first case $\Mu=\H_{(0,-1,0^{H-3},1)}(\Lm)$ and in the second case
$\Mu=\H_{(-1,0^{H-3},1,0)}(\Lm)$.
Note that in all cases the only bead of $\tilde\Lm$,
belonging to the same runner as $a$ is $c$, which is not normal.
Therefore there is no $\nu$-good node of residue $\res A$ and
$U(\nu,\mu)\le m_A(\nu,\mu)=\epsilon(\nu,\mu_A)=1$.

Now consider the case $\epsilon(\nu,\mu_A)=0$.
If there is no $\nu$-good node of residue $\res A$,
then $U(\nu,\mu)=m_A(\nu,\mu)=0$.
Therefore we assume that there is a $\nu$-good node $B$
of residue $\res A$.
By the choice of $A$, the partition $\nu_B$ is almost completely splittable.
Applying Lemma~\ref{l:6:1}\ref{l:6:1:c:2},
we get $\nu_B\not\ged\mu_A$, whence $(\nu_B,\mu_A)\in X$.
Hence $U(\nu,\mu)\le m_A(\nu,\mu)=U(\nu_B,\mu_A)$.
Therefore in the sequel we will consider the case $U(\nu_B,\mu_A)>0$.

Let $b$ be the bead of $\tilde\Lm$ corresponding to $B$.
If $\chi(\lm)=p$ then $b$ is different from the greatest bead of
$\tilde\Lm$ (i.e., $r(B)>1$),
since otherwise the greatest bead of $\tilde\Lm_b$
equal to $b-1$ is not movable up, which contradicts
the fact that $\nu_B$ is almost completely splittable.
Hence $\nu_B=\widetilde{\lm_D}$,
where $D$ is a $\lm$-good node
(see the definition of $\H_\epsilon$ or~\cite[Lemma~8]{Shchigolev11}).

We put $\bar H=h(\lm_D)$.
By the inductive hypothesis, we get
$\mu_A=\H_\epsilon(\lm_D)$, where $\epsilon$ is the sequence
described in cases~\ref{t:upperest:1:c:1}-\ref{t:upperest:1:c:4}
of the current theorem, where $H$ is replaced by $\bar H$.
Since $\sum_{i=1}^{\bar H}\epsilon_i=0$, we have
$\Mu_a=\H_\epsilon(\Lm_d)$, where $d$ is a bead of $\Lm$
corresponding to $D$.

Recall that $a$ and $d$ are cogood spaces of $\Mu_a$ and $\Lm_d$ respectively,
since $a$ and $d$ are cogood beads of $\Mu$ and $\Lm$ respectively.
Since $\epsilon\ne(0^{\bar H})$, we have $\epsilon_j<0$ for some
$1\le j\le\bar H$.
Hence $e\le\b^{\Lm_d}(j)-p$, where
$e$ is the smallest space of $\Lm_d$.
Therefore $d\ne e$,
as otherwise we would get $\b^\Lm=\b^{\Lm_d}$, $\b_\Lm=e$,
$\b^\Lm-\b_\Lm\ge\b^{\Lm_d}(j)-e\ge p$,
which contradicts the fact that $\Lm$ is completely splittable.
Hence $H=\bar H$.

We have $d=\b^{\Lm_d}(k)+1<\b_{\Lm_d}+p$,
as $\b_{\Lm_d}+p$ is not a conormal space of $\Lm_d$.
Obviously the only cogood space of $\Mu_a$,
which equals $\H_\epsilon(\Lm_d)$
from the same runner as $d$ is $\b^{\Lm_d}(k)+p\epsilon_k+1$.
Hence $a=\b^{\Lm_d}(k)+p\epsilon_k+1$, $\Mu=\H_\epsilon(\Lm)$
and $U(\nu,\mu)\le1$ in cases~\ref{t:upperest:1:c:1}-\ref{t:upperest:1:c:3}.

Case 3: $\mu\ged\nu$ and $(\nu,\mu)$ is minimal.
The desired assertion follows from the definitions of $U$ and $\zeta$.
\enddok

The general form of the sequences $\epsilon(H,x,i)$
is quite complicated.
However we have the following assertion.

\begin{lemma}\label{l:upperest:2}
Let $H$ be an integer such that $\tfrac{p+3}2\le H<p$.
Then $\pi(H,x,i)$ is satisfied if and only if $x=QH+1$, $Q\in\Z$, $Q\ge1$ and
$i=H-1$.
In this case
$\epsilon(H,x,i)=(-Q-1,{\bar q}^{H-1-r},(\bar q+1)^r)$,
where $\bar q=\quo(Q+1,H-1)$ and $r=\rem(Q+1,H-1)$.
\end{lemma}
\Dok
Assume $\pi(H,x,i)$ is satisfied.
We shall use the notation of Definition~\ref{o:auxest:1}.
We put
$$
\delta=H-Q-1+h\left[\tfrac{x-Q-1}{p-H+1}\right]+(R-1)\left[\tfrac{x-Q-2}{p-H+1}\right]-
m\left[\tfrac{p-h+x-1}{p-m}\right].
$$

Suppose $R>2$.
Since $m\le h$, we have
$$
\delta\ge
h\left(\left[\tfrac{x-Q-1}{p-H+1}\right]-\left[\tfrac{x-1}{p-H+2}\right]\right)+
\left(\left[\tfrac{x-Q-2}{p-H+1}+1\right]-\left[\tfrac xH\right]\right)+R-2.
$$
To obtain a contradiction with condition~\ref{pidef:2}
of Definition~\ref{o:auxest:1}, it suffices to prove that
both differences in the
outer brackets are nonnegative.
This follows from the inequalities
$$
\begin{array}{l}
\tfrac{x-Q-1}{p-H+1}-\tfrac{x-1}{p-H+2}=
\tfrac{x-1-Q(p-H+2)}{(p-H+1)(p-H+2)}=
\tfrac{Q(2H-(p+2))+R-1}{(p-H+1)(p-H+2)}>0,\\[12pt]
\tfrac{x-Q-2}{p-H+1}+1-\tfrac xH=
\tfrac{(2H-(p+2))x+H(p-1-H)+R}{H(p-H+1)}>0.
\end{array}
$$

Now suppose $R=2$.
Then $Q\ge1$.
We have
$$
\delta\ge
h\left(\left[\tfrac{x-Q-1}{p-H+1}\right]-\left[\tfrac{x-1}{p-H+2}\right]\right)+
\left(\left[\tfrac{x-Q-2}{p-H+1}\right]-\left[\tfrac xH\right]\right)+R-1.
$$
Nonnegativity of the difference in the first pair of the brackets is shown
just as above.
We have
$$
\tfrac{x-Q-2}{p-H+1}-\tfrac xH=\tfrac{(2H-(p+3))x+H(Q-2)+4}{H(p-H+1)}.
$$
The last expression and thus the difference in the second pair of the brackets
is nonnegative if $Q\ge2$.
If $Q=1$ then $x=H+2$ and we have
$$
\left[\tfrac{x-Q-2}{p-H+1}\right]-\left[\tfrac xH\right]=
\left[\tfrac{H-1}{p-H+1}\right]-\left[\tfrac {H+2}H\right]=
\left[\tfrac{H-1}{p-H+1}\right]-1\ge0.
$$
We have a contradiction with condition~\ref{pidef:2}
of Definition~\ref{o:auxest:1}.

Thus we have proved that $R=1$.
Hence $x=QH+1$, $Q\ge1$ and $m=i$.
Suppose $i\ne H-1$.
Then $i\le H-2$.
We have
$$
\delta{\ge}
H-Q-1+h\left[\tfrac{x-Q-1}{p-H+1}\right]-
(H-2)\left[\tfrac{p-H+x}{p-H+2}\right]\!=
h\( \left[\tfrac{x-Q-1}{p-H+1}\right]\!-\!\left[\tfrac{x-2}{p-H+2}\right] \)+
\(\left[\tfrac{x-2}{p-H+2}\right]\!-\!\left[\tfrac xH\right]\)+1.
$$
The difference in the first pair of the outer brackets of the right hand side
is nonnegative, as
$$
\tfrac{x-Q-1}{p-H+1}-\tfrac{x-2}{p-H+2}=
\tfrac{x-2-(Q-1)(p-H+2)}{(p-H+1)(p-H+2)}=
\tfrac{(Q-1)(2H-(p+2))+H-1}{(p-H+1)(p-H+2)}>0.
$$
We have
$$
\tfrac{x-2}{p-H+2}-\tfrac xH=\tfrac{(2H-(p+2))x-2H}{H(p-H+2)}=
\tfrac{(2H-(p+3))x+H(Q-2)+1}{H(p-H+2)}.
$$
The last expression and thus the difference in the second pair of
the brackets is nonnegative if $Q\ge2$.
If $Q=1$ then $x=H+1$ and we have
$$
\left[\tfrac{x-2}{p-H+2}\right]-\left[\tfrac xH\right]=
\left[\tfrac{H-1}{p-H+2}\right]-\left[\tfrac {H+1}H\right]=
\left[\tfrac{H-1}{p-H+2}\right]-1\ge0.
$$
We have a contradiction with condition~\ref{pidef:2}
of Definition~\ref{o:auxest:1}.

Finally for $x=QH+1$ and $i=H-1$, where $Q$ is any positive integer,
we have
$$
\delta=h\(\left[\tfrac{x-1-Q}{p-H+1}\right]-\left[\tfrac{x-1}{p-H+1}\right]\)-Q<0,
$$
whence it follows that $\pi(H,x,i)$ is satisfied.
Now it is clear that $\epsilon(H,x,i)$ is given exactly
by the suggested formula.

\enddok

\section{Applications to branching rules}\label{modstr}

\subsection{Preliminary facts}\label{prelim}
In the following proposition, which follows directly
from~\cite[Theorems~E, $\rm E'$]{Kleshchev_tf3},
$\Res_\alpha0$ is understood as $0$.
\begin{utv}\label{utv:prelim:1}
Let $\lm$ be a $p$-regular partition and $\alpha\in\Z_p$.
\begin{itemize}
\item
If there is no $\lm$-normal {\rm(}$\lm$-conormal{\rm)} node
of residue $\alpha$, then $\Res_\alpha D^\lm=0$
{\rm(}$\Ind^\alpha D^\lm=0${\rm)}.
\item
If there is exactly one $\lm$-normal {\rm(}$\lm$-conormal{\rm)} node
$A$ of residue $\alpha$,
then $\Res_\alpha D^\lm\cong D^{\lm_A}$
{\rm(}$\Ind^\alpha D^\lm\cong D^{\lm^A}${\rm)}.
\end{itemize}
\end{utv}

\begin{lemma}\label{l:prelim:1}
Let $M$ be a module such that $\head M\cong\soc M$ and
$N_1,\ldots,N_k$ be mutually nonisomorphic simple modules
such that $[M:N_i]=1$, $i=1$, \ldots, $k$.
Then either $M\cong\bigoplus_{i=1}^kN_i$ or
there is a simple module $N$ nonisomorphic to any of
$N_1$, \ldots, $N_k$ such that $\Hom(M,N)\ne0$.
\end{lemma}
\Dok
Suppose there is no such a module $N$.
Prove by induction on $n=0$, \ldots, $k$ that there exists a
subset $S_n\subset\{1$, \ldots, $k\}$ of cardinality $n$ such that
$\bigoplus_{i\in S_n}N_i$ is isomorphic to a submodule of $M$.
The set $S_0=\ems$ corresponds to the case $n=0$.
Now let $0<n<k$ and $\iota:\bigoplus_{i\in S_n}N_i\to M$
be an embedding of modules.
Since $|S_n|<k$, we have $\im\iota\ne M$ and there is
a maximal submodule $M_0$ of $M$ containing $\im\iota$.
By our assumption, $M/M_0\cong N_j$ for some
$j\in\{1$, \ldots, $k\}\setminus S_n$.
Since $\head M\cong\soc M$, we can put $S_{n+1}=S_n\cup\{j\}$.

We have $S_k=\{1$, \ldots, $k\}$.
By our assumption from the beginning of the proof, we get
$M\cong\bigoplus_{i=1}^kN_i$.
\enddok

\begin{lemma}\label{l:prelim:2}
Let $\lm$ be a partition of height less than $p$
and $B$ be a $\lm$-addable node such that
$\lm^B$ is $p$-singular.
Then $\lm=(1^{p-1})$ and $B=(p,1)$.
\end{lemma}
\Dok
Since $h(\lm^B)\le p$ and $\lm^B$ is $p$-singular,
we have $\lm^B=(k^p)$.
Hence $k=(\lm^B)_p=1$ and $B=(p,1)$.
\enddok

\subsection{Inducing completely splittable modules}\label{indcs}
In what follows, the image of an $R$-module $M$
in the Grothendieck group of $R$ is denoted by $[M]$.

\begin{teo}\label{t:indcs:1}
Let $\lm$ be a completely splittable partition of $n$ and $\alpha\in\Z_p$.
Suppose there are more than one $\lm$-conormal nodes of residue $\alpha$.
Then there are exactly two such nodes.
Denote them by $A$ and $B$, where $A$ is above $B$.
Then {\rm(}in the Grothendieck group of $K\Sigma_{n+1}${\rm)}
      $$
      [\Ind^\alpha D^\lm]=
      \left\{
      \begin{array}{ll}
      2[D^{\lm^A}]+[D^{\lm^B}]&\!\mbox{if }h_{1,1}(\lm)\ne p-1\mbox{ or }p>2,\lm=(p-1);\\[8pt]
      2[D^{\lm^A}]+[D^{\lm^B}]+[D^{\widetilde{\lm^A}}]&\!\mbox{if }h_{1,1}(\lm)=p-1\mbox{ and }\lm\ne(1^{p-1}),(p-1);\\[8pt]
      2[D^{(2,1^{p-2})}]+[D^{(3,1^{p-3})}]&\!\mbox{if }p>2\mbox{ and }\lm=(1^{p-1});\\[8pt]
      2[D^{(2)}]&\!\mbox{if }p=2\mbox{ and }\lm=(1).
      \end{array}
      \right.
      $$
\end{teo}
\Dok
There are exactly two such nodes, since all proper beads of
any abacus of $\lm$ belong to different runners.
We put $\nu=\lm^A$.
Clearly, $\nu$ is completely splittable,
$A$ is a $\nu$-good node and $h(\lm)=h(\nu)$.
Suppose the assertion of the current theorem does not hold.

First consider the case where $h_{1,1}(\lm)\ne p-1$ or $p>2$, $\lm=(p-1)$.
Then $\lm^B$ is $p$-regular by Lemma~\ref{l:prelim:2} and
$[\Ind^\alpha D^\lm:D^{\lm^B}]=1$.
Applying Lemma~\ref{l:prelim:1} to the following parameters
$$
\hbox to \textwidth
{
\hfil
$M:=\rad\Ind^\alpha D^\lm/\soc\Ind^\alpha D^\lm$,\hfil
$k:=1$,\hfil
$N_1:=D^{\lm^B}$,
\hfil
}
$$
and taking into account $[\Ind^\alpha D^\lm:D^{\lm^A}]=2$,
we get that there is a module $D^\gamma$ isomorphic neither to
$D^{\lm^B}$ nor $D^{\lm^A}$ such that
$\Hom_{\Sigma_{n+1}}(\rad\Ind^\alpha D^\lm,D^\gamma)\ne0$.

By Lemma~\ref{l:spechtex:2.5}, we have
$
0<[\Ind^\alpha D^\lm:D^\gamma]\le2[S^{\lm^A}:D^\gamma]
$
and thus $\nu=\lm^A\led\gamma$.
Hence $\lm\ne(p-1)$ (and therefore $h_{1,1}(\lm)\ne p-1$),
$h(\nu)>1$ and $p>2$.
Since $\head\Ind^\alpha D^\lm\cong D^\nu$, we have
$\Ext_{\Sigma_{n+1}}(D^\nu,D^\gamma)\ne0$.
Theorem~\ref{t:cs:1} implies $\gamma=\tilde\nu$.
We have $h(\lm)=h(\nu)>1$ and $p\le h_{1,1}(\nu)\le h_{1,1}(\lm)+1$
by Lemma~\ref{l:cs:3}.
Taking into account $h_{1,1}(\lm)\ne p-1$, we get $h_{1,1}(\lm)\ge p$ and
thus $\lm$ is a big partition.
Therefore $\tilde\nu=\tilde\lm^C$, where $C$ is a
$\tilde\nu$-good node of residue $\alpha$.
Since $\Ext_{\Sigma_{n+1}}(S^\nu,D^{\tilde\nu})=0$
by Proposition~\ref{u:notation:0.25}, we have by Lemma~\ref{l:radind:2} that
$$
1{\le}\dim\Hom_{\Sigma_{n+1}}(\rad\Ind^\alpha D^\lm,D^\gamma){=}
\dim\Hom_{\Sigma_{n+1}}(\rad\Ind^\alpha D^{\nu_A},D^{\tilde\nu}){\le}
\varepsilon(\nu,{\tilde\nu}_C){=}\varepsilon(\lm^A,\tilde\lm).
$$
Hence $\lm^A=\tilde\lm^D$.
This is possible only if $r(A)=1$ and $\chi(\lm)=p$.
This contradicts the fact that $A$ is a $\lm$-conormal node.

Now consider the case where $h_{1,1}(\lm)=p-1$ and $\lm\ne(1^{p-1}),(p-1)$.
We have $p>2$, $n\ge2$ and $A=(1,\lm_1)$.
We put $A^t=(\lm_1,1)$.
We have $(\lm^t)^{A^t}=(\lm^A)^t$.
Since $\lm$ is a $p$-core,
by~\cite[Lemma~5.2]{Mullineux1}, we get $m(\lm)=\lm^t$.
By~\cite[Lemma~11]{Shchigolev11}, we get $m((\lm^A)^t)=\widetilde{\lm^A}$.
By~\cite[Theorem E(iv)]{Kleshchev_tf3} and Proposition~\ref{u:est:1},
we get
$$
\begin{array}{l}
1=[D^{\lm^t}\uparrow^{\Sigma_{n+1}}:D^{(\lm^t)^{A^t}}]=
[D^{\lm^t}\uparrow^{\Sigma_{n+1}}\otimes\sgn_{n+1}:D^{(\lm^A)^t}\otimes\sgn_{n+1}]=\\[8pt]
[D^{m(\lm^t)}\uparrow^{\Sigma_{n+1}}:D^{m((\lm^A)^t)}]=
[D^\lm\uparrow^{\Sigma_{n+1}}:D^{\widetilde{\lm^A}}]=
[\Ind^\alpha D^\lm:D^{\widetilde{\lm^A}}].
\end{array}
$$
In the case under consideration $\lm^B$ is $p$-regular and
$[\Ind^\alpha D^\lm:D^{\lm^B}]=1$.
Applying Lemma~\ref{l:prelim:1} to the following parameters
$$
\hbox to \textwidth
{
\hfil
$M:=\rad\Ind^\alpha D^\lm/\soc\Ind^\alpha D^\lm$,\hfil
$k:=2$,\hfil
$N_1:=D^{\lm^B}$,\hfil
$N_2:=D^{\widetilde{\lm^A}}$
\hfil
}
$$
and taking into account $[\Ind^\alpha D^\lm:D^{\lm^A}]=2$,
we get that there is a module $D^\gamma$ isomorphic to
none of $D^{\lm^B}$, $D^{\widetilde{\lm^A}}$, $D^{\lm^A}$,
such that $\Hom_{\Sigma_{n+1}}(\rad\Ind^\alpha D^\lm,D^\gamma)\ne0$.
Since $\head\Ind^\alpha D^\lm\cong D^{\lm^A}$, we have
$\Ext_{\Sigma_{n+1}}(D^{\lm^A},D^\gamma)\ne0$.
Similarly to the previous case Lemma~\ref{l:spechtex:2.5} implies
$\lm^A\led\gamma$.
Hence by Theorem~\ref{t:cs:1},
we get a contradiction $\gamma=\widetilde{\lm^A}$.

If $p>2$ then multiplying $[\Ind^{\overline{-1}}D^{(p-1)}]\cong2[D^{(p)}]+[D^{(p-1,1)}]$
by $\sgn_p$ and applying Proposition~\ref{u:est:1}, we get
$[\Ind^{\bar1}D^{(1^{p-1})}]\cong2[D^{(2,1^{p-2})}]+[D^{(3,1^{p-3})}]$.

Finally the formula $[\Ind^{\bar1}D^{(1)}]=2[D^{(2)}]$ for $p=2$
can be checked by dimension comparison.
\enddok

\begin{teo}\label{t:indcs:2}
Let $\lm$ be a completely splittable partition of $n$ different
from $(1^{p-1})$ and $\alpha\in\Z_p$.
Denote by $B$ the bottom $\lm$-addable node
(i.e. from the first column).
Suppose there are more than one $\lm^B$-normal nodes
of residue $\alpha$.
Then except $B$ there is only one such node $A$.
We have $\Res_\alpha D^{\lm^B}\cong\Ind^\alpha D^{\lm_A}$.
\end{teo}
\Dok
Let $B_1$, \ldots, $B_k$ be all $\lm$-conormal nodes
different from $B$.
Since there are more than one $\lm^B$-normal nodes
of residue $\alpha$, the residues $\res B$, $\res B_1$, \ldots, $\res B_k$
are mutually distinct.
Hence
$D^\lm\uparrow_{\Sigma_{n+1}}\cong D^{\lm^B}\oplus D^{\lm^{B_1}}\oplus\cdots\oplus D^{\lm^{B_k}}$.
Now it is clear that the only nonsimple indecomposable summand of
$D^\lm\uparrow_{\Sigma_{n+1}}\downarrow_{\Sigma_n}$ is
$\Res_\alpha D^{\lm^B}$.

Let $A_1$, \ldots, $A_l$ be all the $\lm$-normal nodes.
The residues $\res A_1$, \ldots, $\res A_l$ are mutually distinct and
$D^\lm\downarrow_{\Sigma_{n-1}}\cong D^{\lm_{A_1}}\oplus\cdots\oplus D^{\lm_{A_l}}$.
The only nonsimple indecomposable summand of
$D^\lm\downarrow_{\Sigma_{n-1}}\uparrow_{\Sigma_n}$ is
$\Ind^\alpha D^{\lm_A}$.

By the subgroup theorem~\cite[(44.2)]{Curtis_Reiner}
applied to $G:=\Sigma_{n+1}$, $R=S:=\Sigma_n$ and $L:=D^\lm$,
we have $D^\lm\uparrow_{\Sigma_{n+1}}\downarrow_{\Sigma_n}\cong%
D^\lm\downarrow_{\Sigma_{n-1}}\uparrow_{\Sigma_n}\oplus D^\lm$.
By the Krull-Schmidt theorem, we have
$\Res_\alpha D^{\lm^B}\cong\Ind^\alpha D^{\lm_A}$.
\enddok

\subsection{Inducing some almost completely splittable modules}\label{indacs}

\begin{teo}\label{t:indacs:2}
Let $\lm$ be a big partition of $n$ having height $H\ge\tfrac{p+3}2$ and
$\alpha\in\Z_p$ such that
$h_{2,1}(\lm)\ne p-1$ and
the condition $\overline{\lm_1}=-\overline{h(\lm)}=\alpha$ does not hold.
Suppose there are more than one $\tilde\lm$-conormal nodes
of residue $\alpha$.
Then there are exactly two such nodes.
Denote them by $A$ and $B$, where $A$ is above $B$.
We have
$
[\Ind^\alpha D^{\tilde\lm}]=2[D^{\tilde\lm^A}]+[D^{\tilde\lm^B}]
$
{\rm(}in the Grothendieck group of $K\Sigma_{n+1}${\rm)}.
\end{teo}
\Dok
Arguing as in the proof of Theorem~\ref{t:indcs:1},
we obtain that there exactly two such nodes.
Let $\Lm$ be an abacus of $\lm$ and $c$ be its minimal space.
Denote by $a$ and $b$ the spaces of $\tilde\Lm$ corresponding to
$A$ and $B$ respectively.
The case $h(\tilde\lm)<h(\lm)$ is impossible, as we would have $c=\b^\Lm-p$,
which contradicts the existence of more than one initial space of $\tilde\Lm$
in the same runner.
Therefore $h(\tilde\lm)=h(\lm)$ and $c=b$.
Hence $\tilde\lm$ is not completely splittable and
in particular $\tilde\lm\ne(1^{p-1})$.
By Lemma~\ref{l:prelim:2}, $\tilde\lm^B$ is $p$-regular.

Since $\res A=\res B$ and
the condition $\overline{\lm_1}=-\overline{h(\lm)}=\alpha$ is violated,
there exists $i=2$, \ldots, $H$ such that the runner
containing $\b^\Lm(i)+1$ contains the space $c$ and
no proper bead of $\Lm$.
Hence $\tilde\lm^A=\widetilde{\lm^D}$,
where $D$ is a $\lm$-cogood node such that $h(\lm^D)=H$.
Therefore $\lm^D$ is big.

Suppose that the assertion of the current theorem does not hold.
We put $\nu=\tilde\lm^A$ and denote by $d$ the space of $\Lm$
corresponding to $D$.
Applying Lemma~\ref{l:prelim:1} to the parameters
$$
\hbox to \textwidth
{
\hfil
$M:=\rad\Ind^\alpha D^{\tilde\lm}/\soc\Ind^\alpha D^{\tilde\lm}$,\hfil
$k:=1$,\hfil
$N_1:=D^{\tilde\lm^B}$
\hfil
}
$$
and taking into account $[\Ind^\alpha D^{\tilde\lm}:D^{\tilde\lm^A}]=2$,
we get that there is a module $D^\gamma$ isomorphic neither to
$D^{\tilde\lm^B}$ nor $D^{\tilde\lm^A}$ such that
$\Hom_{\Sigma_{n+1}}(\rad\Ind^\alpha D^{\tilde\lm},D^\gamma)\ne0$.
By Lemma~\ref{l:spechtex:2.5}, we get
$
0<[\Ind^\alpha D^{\tilde\lm}:D^\gamma]\le2[S^{\tilde\lm^A}:D^\gamma]
$
and thus $\nu=\tilde\lm^A\led\gamma$.
Since $\head\Ind^\alpha D^{\tilde\lm}\cong D^\nu$, we have
$\Ext_{\Sigma_{n+1}}(D^\nu,D^\gamma)\ne0$.
Theorem~\ref{t:upperest:1} and Lemma~\ref{l:upperest:2}
imply $\gamma=\H_\epsilon(\lm^D)$ for
$\epsilon$ equal to one of the following sequences:
$(0,-1,0^{H-3},1)$; $(-1,0^{H-3},1,0)$; $(-1,-1,0^{H-4},1,1)$;
$(-Q-1,{\bar q}^{H-1-r},(\bar q+1)^r)$,
where $Q\ge1$, $\bar q=\quo(Q+1,H-1)$ and $r=\rem(Q+1,H-1)$.

If we suppose that $h(\gamma)<H$, then,
taking into account the exact form of
possible values of $\epsilon$ mentioned above,
we get $c=\b^{\Lm^d}-p$ or $c=\b^{\Lm^d}(2)-p$.
However the former condition does not hold as $i>1$ and
the latter does not hold as $h_{2,1}(\lm)\ne p-1$.
Hence $h(\gamma)=H$,
there exists a unique $\gamma$-normal node $E$ of residue $\alpha$
and $\H_\epsilon(\lm)=\gamma_E$.

Since $\Ext_{\Sigma_{n+1}}(S^\nu,D^\gamma)=0$
by Proposition~\ref{u:notation:0.25}, we have by Lemma~\ref{l:radind:2} that
$$
1{\le}\dim\Hom_{\Sigma_{n+1}}(\rad\Ind^\alpha D^{\tilde\lm},D^\gamma){=}
\dim\Hom_{\Sigma_{n+1}}(\rad\Ind^\alpha D^{\nu_A},D^\gamma){\le}
\varepsilon(\tilde\lm^A,\gamma_E).
$$
Hence $\tilde\lm^A=(\gamma_E)^F$, where $\res F=\alpha$.
Clearly, either $F=E$ or $h((\gamma_E)^F)>H$.
The former case is impossible as $\tilde\lm^A\led\gamma$
and the latter is impossible as $h(\tilde\lm^A)=H$.
\enddok

{\bf Remark.} The reader can by oneself formulate and prove
the analog of Theorem~\ref{t:indcs:2}
for almost completely splittable partitions.

\subsection{Conjectures}\label{hypothesis}
The following conjectures are based on
Theorem~\ref{t:upperest:1} and calculations
within the known decomposition matrices.

\begin{gip}\label{gip:hypothesis:1}
Let $\lm$ be a big partition of $n$ having height $H$ and
$\alpha\in\Z_p$ such that
$\overline{\lm_1}=-\overline{h(\lm)}=\alpha$ does not hold.
Suppose there are more than one $\tilde\lm$-conormal nodes
of residue $\alpha$.
Then there are exactly two such nodes.
Denote them by $A$ and $B$, where $A$ is above $B$.
We have
$$
\begin{array}{l}
\!
\!
\!
[\Ind^\alpha D^{\tilde\lm}]=\\[12pt]
\;
\;
\left\{\!\!
\begin{array}{ll}
2[D^{\tilde\lm^A}]+[D^{\tilde\lm^B}]&\!\!\!\mbox{if }h_{2,1}(\lm){\ne}p-1;\\[8pt]
2[D^{\tilde\lm^A}]+[D^{\tilde\lm^B}]+[D^{\H_{(-2,2)}(\tilde\lm^A)}]&\!\!\!\mbox{if }h_{2,1}(\lm){=}p-1\mbox{ and }H{=}2;\\[8pt]
2[D^{\tilde\lm^A}]+[D^{\tilde\lm^B}]+[D^{\H_{(0,-1,1)}(\tilde\lm^A)}]        +[D^{\H_{(1,-1,0)}(\tilde\lm^A)}]&\!\!\!\mbox{if }h_{2,1}(\lm){=}p-1\mbox{ and }H{=}3;\\[8pt]
2[D^{\tilde\lm^A}]+[D^{\tilde\lm^B}]+[D^{\H_{(0,-1,0^{H-3},1)}(\tilde\lm^A)}]+[D^{\H_{(0,-1,0^{H-4},1,0)}(\tilde\lm^A)}]&\!\!\!\mbox{if }h_{2,1}(\lm){=}p-1\mbox{ and }H{>}3.
\end{array}
\!
\!
\!
\!
\!
\!
\!
\!\!\!
\right.
\end{array}
$$
\end{gip}
{\it Example.}
Let $p=7$ and $\lm=(5,4,2,2)$.
Then $\tilde\lm=(6,4,2,1)$ and
$[\Ind^3D^{(6,4,2,1)}]=2[D^{(6,5,2,1)}]+[D^{(6,4,2,1,1)}]+[D^{(6,6,2)}]+[D^{(7,7)}]$.
We have $\lm_1=5\ne-4=-h(\lm)\pmod7$, $h_{2,1}(\lm)=6$, $A=(2,5)$, $B=(5,1)$,
$\tilde\lm^A=(6,5,2,1)$, $\tilde\lm^B=(6,4,2,1,1)$,
$\H_{(0,-1,0,1)}(\tilde\lm^A)=(6,6,2)$,
$\H_{(0,-1,1,0)}(\tilde\lm^A)=(7,7)$.

\begin{gip}\label{gip:hypothesis:2}
Let $\lm$ be a completely splittable partition of height $H$
such that $\chi(\lm)=p$ and let $\alpha\in\Z_p$.
Denote by $A$ the top $\lm$-addable node
(i.e. from the first row).
Suppose there are more than one $\lm^A$-normal nodes of residue $\alpha$.
Then except $A$ the only such node is the bottom $\lm$-removable node $B$.
We have
$$
[\Res_\alpha D^{\lm^A}]=\left\{
\begin{array}{ll}
2[D^{\tilde\lm}]+[D^\lm]+[D^{\H_{(-1,0^{H-3},1,0)}(\lm)}]+
x[D^{   \H_{(0,-1,0^{H-3},1)}(\lm)   }]&\mbox{if }H>2;\\
2[D^{\tilde\lm}]+[D^\lm]+
x[D^{   \H_{(1,-1)}(\lm)   }]&\mbox{if }H=2,
\end{array}
\right.
$$
where $x=[h_{2,1}(\lm)\ge p]$.
\end{gip}
{\it Example.}
Let $p=5$ and $\lm=(5,5,3)$.
Then $A=(1,6)$, $\lm^A=(6,5,3)$ and
$[\Res_0D^{(6,5,3)}]=2[D^{(6,5,2)}]+[D^{(5,5,3)}]+%
 [D^{(9,2,2)}]+[D^{(6,6,1)}]$.
We have $\tilde\lm=(6,5,2)$, $\H_{(-1,1,0)}(\lm)=(9,2,2)$,
$\H_{(0,-1,1)}(\lm)=(6,6,1)$ and $h_{2,1}(\lm)=6$.

\begin{gip}\label{gip:hypothesis:3}
If $\lm$ is a completely splittable partition of height $3$
such that $h_{1,1}(\lm)=2p-1$, then
$[\Ind^{-\bar3}D^{\tilde\lm}]=
2[D^{\tilde\lm^A}]+[D^{\tilde\lm^B}]+
[D^{\H_{(0,1,-1)}(\tilde\lm^A)}]$,
where $A=(3,\tilde\lm_3+1)$ and $B=(4,1)$.
\end{gip}
{\it Example.}
Let $p=5$ and $\lm=(7,6,6)$.
Then $\tilde\lm=(9,6,4)$ and
$[\Ind^2D^{(9,6,4)}]=2[D(9,6,5)]+[D(9,6,4,1)]+[D(10,10)]$.
We have $h_{1,1}(\lm)=9$, $\tilde\lm^A=(9,6,5)$,
$\tilde\lm^B=(9,6,4,1)$ and $\H_{(0,1,-1)}(\tilde\lm^A)=(10,10)$.
Moreover with the help of the known decomposition matrices and
the subgroup theorem~\cite[(44.2)]{Curtis_Reiner}, it is easy to
verify Conjecture~\ref{gip:hypothesis:3} for $p=5$.

\end{document}